\documentclass[leqno]{article}
\usepackage{amssymb, amsmath, amsfonts}

\newtheorem{theorem}{Theorem}[section]
\newtheorem{corollary}[theorem]{Corollary}
\newtheorem{lemma}[theorem]{Lemma}
\newtheorem{proposition}[theorem]{Proposition}

\numberwithin{equation}{section}

\def\sqr#1#2{{\vcenter{\vbox{\hrule height.#2pt
    \hbox{\vrule width.#2pt height#1pt \kern#1pt
    \vrule width.#2pt}
    \hrule height.#2pt}}}}

\def\ecarre{\; \mathchoice\sqr56\sqr56\sqr{4.3}5\sqr{3.7}5}

\def\Ga{\Gamma}
\def\ga{\gamma}
\def\al{\alpha}
\def\De{\Delta}
\def\de{\delta}
\def\La{\Lambda}
\def\la{\lambda}
\def\eps{\varepsilon}
\def\ka{\kappa}

\def\om{\omega}
\def\Om{\Omega}

\def\pa{\partial}

\def\Dom{\hbox{Dom}}

\def\tB{\tilde B}

\def\ta{\tilde a}

\font\bbb=msbm10

\def\R{\hbox{{\bbb R}}}

\def\cL{{\cal L}}

\def\Dom{\, {\rm Dom}\, }
\def\capa{\, {\rm cap}\, }
\def\mes{\,{\rm mes}\, }
\def\Lip{\, {\rm Lip}\, }
\def\Liploc{\,{\rm Lip}_{\textit{loc}}\,}

\def\Tr{\, {\rm Tr}\, }

\def\curl{\hbox{curl}\;}
\def\supp{\hbox{supp}\;}
\def\lBr{\langle B\rangle}
\def\lBxr{\langle B(x)\rangle}
\def\lByr{\langle B(y)\rangle}
\def\Veff{V_{\textit{eff}}}
\def\Veffde{V_{\textit{eff}}^{(\de)}}
\def\Veffdeeps{V_{\textit{eff}}^{(\de,\eps)}}

\def\ni{\noindent}

\def\ms{\medskip}
\def\bs{\bigskip}

\begin{document}
\title{Discreteness of spectrum for the magnetic Schr\"odinger
operators. I.}
\author{{\bf Vladimir Kondratiev\thanks{Research partially supported by 
Oberwolfach Forschungsinstitut f\"ur Mathematik }}
\smallskip\\
Department of Mechanics and Mathematics\\
Moscow State University\\
Vorobievy Gory, Moscow, 119899,
Russia\medskip\\
E-mail: kondrat@vnmok.math.msu.su
\bigskip\\
{\bf Mikhail Shubin
\thanks{Research partially supported by NSF grant DMS-9706038 and
Oberwolfach Forschungsinstitut f\"ur Mathematik } }
\\
Department of Mathematics\\
Northeastern University\\
Boston, MA 02115,
USA
\medskip\\
E-mail: shubin@neu.edu
}
\date{}
\maketitle

\begin{abstract}

We consider a magnetic Schr\"odinger  operator $H$
in $\R^n$ or on a Riemannian manifold $M$ of bounded
geometry.  Sufficient conditions for the  spectrum
of $H$ to be discrete are given in terms of behavior
at infinity for some effective potentials $V_{eff}$
which are expressed through electric and magnetic fields.
These conditions 
can be formulated in the form $V_{eff}(x)\to +\infty$
as $x\to\infty$. They generalize the classical
result by K.Friedrichs (1934), and include earlier results
of J.~Avron, I.~Herbst and B.~Simon (1978), 
A.~Dufresnoy (1983) and A.~Iwatsuka (1990) which
were obtained in the absence of an electric field. 
More precise sufficient conditions
can be formulated in terms of the Wiener capacity
and extend earlier work by A.M.~Molchanov (1953) and
V.~Kondrat'ev and M.~Shubin (1999)
who considered the case of the operator without a magnetic field.
These conditions become necessary and sufficient in case there is
no magnetic field and the electric potential is semi-bounded below.

\end{abstract}

\section{Introduction and main results}\label{S:Intro}

\subsection{Notations and preliminaries}\label{SS:notation}

The main object of this paper is a magnetic Schr\"odinger operator
in $\R^n$ and its generalizations. In the simplest case it has the
form

\begin{equation}\label{E:Schroed}
H_{a,V}=\sum_{j=1}^n P_j^2+V,
\end{equation}
where
\begin{equation}\label{E:Pj}
P_j=\frac{1}{i}\frac{\pa}{\pa x^j}+a_j,
\end{equation}
and $a_j=a_j(x)$, $V=V(x)$, $x=(x^1,\dots,x^n)\in\R^n$.
We assume that $a_j$ and $V$ are real-valued functions.

For simplicity we will assume now that $a_j\in C^1(\R^n)$,
$V\in L^\infty_{loc}(\R^n)$ (the later means that $V$ is measurable and
locally bounded). Then $H_{a,V}$ is well defined on $C_c^\infty(\R^n)$
(the set of all complex-valued $C^\infty$ functions with a compact
support in $\R^n$), and it is an unbounded symmetric operator in
$L^2(\R^n)$.

We will always impose (explicitly or implicitly) conditions which insure that
the operator $H_{a,V}$ is essentially self-adjoint. For instance, the
condition $V\ge 0$ is sufficient (see e.g.\ H.~Leinfelder and C.G.~Simader
\cite{Leinfelder-Simader} where this is proved under most general local
regularity conditions on $a_j$ and $V$). But some negative
potentials (even mildly blowing up to $-\infty$ when $x\to\infty$)
will do as well - see e.g.\ T.~Ikebe and T.~Kato \cite{Ikebe-Kato},
A.~Iwatsuka \cite{Iwatsuka2}, M.~Shubin \cite{Shubin3}
for several versions of this fact. For the sake of convenience
of the reader we give in  Section \ref{S:self} 
a very short proof of the fact which is important for us:
the semi-boundedness below for the operator $H_{a,V}$
(on $C_c^\infty(\R^n)$) implies that it is essentially self-adjoint
(this is an extension of the Povzner--Wienholtz theorem --
see \cite{Povzner}, \cite{Wienholtz}, \cite{Glazman}, \cite{Simader78})
and it is proved for the case of operators on any complete 
Riemannian manifold in \cite{Shubin00}.
We will also denote by $H_{a,V}$ the corresponding self-adjoint operator in
$L^2(\R^n)$.

Actually the condition $V\in L^\infty_{loc}(\R^n)$ is not necessary
for our study. For example, it will be sufficient
to have $V\in L^2_{loc}(\R^n)$ and locally semi-bounded below.
Moreover, it is sufficient to have $V\in L^1_{loc}(\R^n)$ and locally semi-bounded
below. In this case we have to impose conditions which guarantee that the 
corresponding quadratic form $h_{a,V}$ is semi-bounded below and consider the operator 
defined by this form. This does not make any  difference in our arguments
because we can work with the quadratic form only. But we prefer 
not to get the reader distracted by unimportant details.

On the other hand working with capacities usually requires $V$ to be locally 
semi-bounded below. So this condition often can not be removed.

We will say that a self-adjoint operator $H$ in
a Hilbert space $\cal{H}$
has a \textit{discrete spectrum} if 
its spectrum consists of isolated eigenvalues of finite
multiplicities. 
It follows that the only accumulation points of
these eigenvalues can be $\pm\infty$. Equivalently we may say that $H$ 
has a compact resolvent.

Our goal will be to provide conditions (mainly sufficient but
sometimes necessary and sufficient) for the discreteness of the
spectrum of $H_{a,V}$. We will abbreviate the discreteness
of the spectrum of $H_{a,V}$ by writing $\sigma=\sigma_d$.
Actually under the conditions we will impose the operator will be
always semi-bounded below, so the only
accumulation point of the eigenvalues will be $+\infty$.

Let us recall first some facts concerning the Schr\"odinger
operator $H_{0,V}=-\De+V$ without magnetic field (i.e.\ the operator
(\ref{E:Schroed}) with $a=0$).

It is a  classical result of K.~Friedrichs \cite{Friedrichs}
(see also  e.g.\ \cite{Reed-Simon}, Theorem XIII.67, or 
\cite{Berezin-Shubin}, Theorem 3.1) that
the condition
\begin{equation}\label{E:V-to-infty}
V(x)\to+\infty\quad\hbox{as}\quad x\to\infty
\end{equation}
implies  $\sigma=\sigma_d$ (for $H_{0,V}$).

Now assume that
\begin{equation}\label{E:semibound}
V(x)\ge -C
\end{equation}
with a constant $C$, i.e.\
the potential $V$ is semi-bounded below. Without loss of generality
we can assume then that $V\ge 0$.
Let us formulate a simple necessary condition for the discreteness
of the spectrum.
Denote by $B(x,r)$ the open ball with the radius $r>0$ and the
center at $x\in\R^n$. Then $\sigma=\sigma_d$ for $H_{0,V}$ implies
that for every fixed $r>0$
\begin{equation}\label{E:av-V-to-infty}
\int_{B(x,r)}V(y)dy\to+\infty \quad\hbox{as}\quad x\to\infty.
\end{equation}
This observation was made in a remarkable paper by A.~Molchanov
\cite{Molchanov} who proved that in case $n=1$
this condition is in fact necessary and sufficient
(assuming the semi-boundedness of the potential $V$).

More importantly,
A.~Molchanov found  a necessary and sufficient condition for
$\sigma=\sigma_d$ to hold, again assuming (\ref{E:semibound}).
This condition is intermediate
between (\ref{E:V-to-infty}) and (\ref{E:av-V-to-infty}). It is
formulated in terms of the Wiener 
capacity which we will denote $\capa$ (see e.g.\ \cite{Edmunds-Evans, Kondrat'ev-Shubin, Mazya8}
for  necessary properties of the
capacity and more details).  In case $n=2$
the capacity of a set $F\subset B(x,r)$ is always taken relative to
a ball $B(x,R)$ of a fixed radius $R>r$.
(Expositions of Molchanov's and more general results can be found
in \cite{Edmunds-Evans, Kondrat'ev-Shubin, Mazya8}.)

A.~Molchanov proved that  $H_{0,V}$ has a discrete spectrum
if and only if there exist $c>0$ and $r_0>0$ such that for any
$r\in (0,r_0)$
\begin{equation}\label{E:M-cond}
\underset{F}{\inf}
\left\{\left.\int_{B(x,r)\setminus F}V(y)dy \right|\;\capa(F)\le
c\cdot\capa(B(x,r))\right\}
\to+\infty \ \hbox{as}\ x\to\infty.
\tag{$M_c$}
\end{equation}
(In this case we will say that the function $V$ satisfies $(M_c)$, or that
the Molchanov condition $(M_c)$ holds for $V$. Later
we will impose this condition on some other functions.)
Note that $(M_c)$ implies $(M_{c'})$ for any $c'<c$. Hence we can
equivalently write that $(M_c)$ is satisfied for all $c\in(0,c_0)$
with a positive $c_0$. In fact A.~Molchanov provides
a particular value of $c$ (e.g.\ $c=2^{-2n-6}$ would do -- see
\cite{Kondrat'ev-Shubin}), though it is by no means precise.

Note also that $\capa(B(x,r))$ can be explicitly calculated.
It equals $c_nr^{n-2}$ if $n\ge 3$. If $n=2$ then
$\capa(B(x,r))$ asymptotically
equals $c_2\left(\log(1/r)\right)^{-1}$ as $r\to 0$.
Hence in the formulation
of the Molchanov condition $(M_c)$ we can replace $\capa(B(x,r))$
by $r^{n-2}$ if $n\ge 3$ and by $\left(\log(1/r)\right)^{-1}$ if
$n=2$.

A simple argument given in \cite{Avron-Herbst-Simon} (see also
Sect.\ref{S:localization} of this paper) shows that if
$H_{0,V}$ has a discrete
spectrum, then the same is true for $H_{a,V}$ whatever the
magnetic potential $a$. Therefore the condition (\ref{E:M-cond})
together with (\ref{E:semibound}) is sufficient for the
discreteness of spectrum of $H_{a,V}$. This means that a magnetic
field can only improve the situation from our point of view.
Papers by J.~Avron, I.~Herbst and B.~Simon
\cite{Avron-Herbst-Simon}, A.~Dufresnoy \cite{Dufresnoy}
and A.~Iwatsuka \cite{Iwatsuka1} provide some quantitative
results which show that even in case $V=0$ the magnetic field 
can make the spectrum discrete. 
In this paper we will improve
the results of the above mentioned papers. In particular we will
add the capacity into the picture, so in many cases our conditions
become necessary and sufficient in case when there is no magnetic
field, i.e.\ when $a=0$. We will also make both electric and magnetic fields
work together to achieve
the discreteness of spectrum. 

Unfortunately we can not provide efficient necessary and
sufficient conditions of the discreteness of the spectrum
when both fields are present (or even if the magnetic field only
is present). 
The conditions which we can give always contain some hypotheses which
are hard to check (unless $a=0$ when they become trivial).
Some of these conditions will be discussed in a future
continuation of this paper.

It is convenient to consider the magnetic potential as
a 1-form $a$ with components $a_j$:
\begin{equation}\label{E:a-form}
a=a_jdx^j,
\end{equation}
where we use the Einstein summation convention (i.e.\ the summation
over all repeated suffices is understood). Now the
\textit{magnetic field} is a 2-form $B$ which is defined as
\begin{equation}\label{E:B-form}
B=da=\frac{\pa a_j}{\pa x^k}dx^k\wedge dx^j=
\frac{1}{2}B_{jk}dx^j\wedge dx^k,
\end{equation}
where $B_{kj}=-B_{jk}$. Obviously
\begin{equation}\label{E:B-form-2}
B=\sum_{j<k}B_{jk}dx^j\wedge dx^k,
\end{equation}
and
\begin{equation}\label{E:Bjk}
B_{jk}=\frac{\pa a_k}{\pa x^j}-\frac{\pa a_j}{\pa x^k},
\end{equation}
so in the standard vector analysis notation $B=\curl a$.
The functions $B_{jk}$ will be called \textit{components}
of the magnetic field $B$.

In case $n=2$ the magnetic field has essentially one
non-trivial component $B_{12}=-B_{21}$ and in this case
we will denote $B=B_{12}$.

We will need a norm of $B$ which is defined as
\begin{equation}\label{E:abs-B}
|B|=\left(\sum_{j<k}|B_{jk}|^2\right)^{1/2}.
\end{equation}

Note that the components of the magnetic field show up in
the commutation relations
\begin{equation}\label{E:commut}
[P_j,P_k]=\frac{1}{i}B_{jk},
\end{equation}
where $[A,B]=AB-BA$ for operators $A,B$ in the same Hilbert space
(at the moment we assume that all operations are performed on
the same domain, e.g.\ $C_c^\infty(\R^n)$ for $P_j$ and $P_k$).
The relation (\ref{E:commut}) allows to apply uncertainty
principle type arguments in investigating the spectrum.

An important fact is the \textit{gauge invariance} of the spectrum
of $H_{a,V}$: this spectrum does not depend of the choice of
the magnetic potential $a$ provided the magnetic field $B$ is
fixed. Namely, if $a,a'$ are two magnetic potentials with
$da=da'=B$, then $\sigma(H_{a,V})=\sigma(H_{a',V})$ for any $V$.
To see this note that by the Poincar\'e Lemma
(see e.g.\ 4.18 in \cite{Warner})
we have $a'=a+d\phi$, where $\phi\in C^1(\R^n)$ is defined up to
an additive constant and can be assumed real-valued. Then the
corresponding operators
$$
P_j^\prime=\frac{1}{i}\frac{\pa}{\pa x^j}+a_j^\prime
$$
are related with $P_j$ by the formulae
$$
P_j^\prime=e^{-i\phi}P_j e^{i\phi}.
$$
Therefore
$$
H_{a',V}=e^{-i\phi}H_{a,V} e^{i\phi},
$$
and the operators $H_{a',V}$ and $H_{a,V}$ are unitarily
equivalent, hence have the same spectra.

For the weakest requirements on the magnetic potential $a$
the gauge invariance was established by H.~Leinfelder \cite{Leinfelder}.

By this reason it is more natural for spectral theory
to formulate the conditions
on $H_{a,V}$ in terms of $B,V$ rather than $a,V$.

Let us assume that we are given a magnetic potential
$a=a_jdx^j$, $a_j\in C^1(\R^n)$.
For a function $u\in C^1(\R^n)$ (or, more generally,
for a locally Lipschitz function) define
\textit{magnetic differential} as

\begin{equation}\label{E:mag-diff}
d_a u=du+iua\in \La^1(\R^n).
\end{equation}
It is also convenient to identify this complex-valued
1-form with the corresponding complex vector field
which is called \textit{magnetic gradient}:

\begin{equation}\label{E:mag-grad}
\nabla_a u=
\left(\frac{\pa u}{\pa x^1}+ia_1u,\dots, \frac{\pa u}{\pa x^n}+ia_nu\right)
=\left(iP_1u,\dots,iP_nu\right).
\end{equation}

We will denote by $|\cdot|$ the usual euclidean norm of vectors
or 1-forms.

\subsection{Localization} \label{SS:localization}

Necessary and sufficient conditions of  discreteness of spectrum
for $H_{a,V}$ can be formulated in terms of bottoms of
Dirichlet or Neumann spectra on balls of a fixed radius
or on cubes of a fixed size. We will call these facts
{\it localization} results. The first result of this kind about
usual Schr\"odinger operators (without magnetic field) is due
to A.~Molchanov \cite{Molchanov} (see also
\cite{Kondrat'ev-Shubin} for a more general theorem
on manifolds). A.~Iwatsuka \cite{Iwatsuka1}
proved a localization theorem
for magnetic Schr\"odinger operators.

The bottoms of Dirichlet and Neumann spectra for
the operator $H_{a,V}$ in an open set $\Om\subset\R^n$
are defined in terms of its quadratic form
which we will denote $h_{a,V}$:

\begin{equation}\label{E:quadr-form}
h_{a,V}(u,u)=\int_\Om (|\nabla_a u|^2+V|u|^2)dx.
\end{equation}
This form is well defined e.g.\ for all $u\in L^2(\Om)$ such that
$P_j u\in L^2(\Om)$,  $j=1,\dots,n$, the derivatives are understood
in the sense of distributions, and $V|u|^2\in L^1(\Om)$. In
particular $h_{a,V}(u,u)$ is well defined for all
$u\in C_c^\infty(\Om)$. Denote by $(u,v)$ the usual scalar product
of $u$ and $v$ in $L^2(\Om)$.

It is easy to check that the following gauge invariance relation
holds:
\begin{equation}\label{E:gauge-quadr}
h_{a+d\phi,V}(u,u)=h_{a,V}(e^{i\phi}u,e^{i\phi}u),
\end{equation}
for any $\phi\in C^1(\R^n)$ and $u$ as above.

Now we can define
\begin{equation}\label{E:lambda}
\la(\Om;H_{a,V})=\underset{u}{\inf} \left\{\frac{h_{a,V}(u,u)}
{(u,u)},\;
u\in C_c^\infty(\Om)\setminus 0\right\}\;,
\end{equation}

\begin{equation}\label{E:mu}
\mu(\Om;H_{a,V})=\underset{u}{\inf} \left\{\frac{h_{a,V}(u,u)}
{(u,u)},\;
u\in (C^\infty(\Om)\setminus\{0\})\cap L^2(\Om),
\int_\Om V|u|^2dx>-\infty\right\}\;,
\end{equation}
i.e.\ $\la(\Om;H_{a,V})$ and $\mu(\Om;H_{a,V})$ are bottoms of
the  Dirichlet and Neumann spectra (of $H_{a,V}$)
respectively, in the usual variational understanding
(see e.g.\ \cite{Courant-Hilbert}, \cite{Kato1}).

The relation \eqref{E:gauge-quadr} obviously implies that
the numbers $\la(\Om;H_{a,V})$ and

\ni
$\mu(\Om;H_{a,V})$
are gauge invariant, i.e.\ they do not change if we replace
$a$ by $a+d\phi$ for any $\phi\in C^1(\R^n)$.

The following
theorem  slightly extends the result of A.~Iwatsuka \cite{Iwatsuka1}
removing the requirement $V\ge 0$ and allowing non-continuous
minorant functions.

\begin{theorem}\label{T:localization}
The following conditions are equivalent:

$(a)$ $H_{a,V}$ is essentially self-adjoint, semi-bounded below
and has a discrete spectrum;

$(b)$ $\la(B(x,r);H_{a,V})\to +\infty$ as $x\to\infty$ for any fixed
$r>0$;

$(c)$ there exists $r>0$ such that
$\la(B(x,r);H_{a,V})\to +\infty$ as $x\to\infty$;

$(d)$ there exists a real valued function $\La\in C(\R^n)$
such that $\La(x)\to +\infty$ as $x\to\infty$
and the operator inequality
\begin{equation}\label{E:H-ge-La}
H_{a,V}\ge\La(x)
\end{equation}
holds in the sense of quadratic forms $($on $C_c^\infty(\R^n)$$)$;

$(e)$ there exists a measurable function $\La:\R^n\to\R$
such that $\La(x)\to +\infty$ as $x\to\infty$ and
(\ref{E:H-ge-La}) holds.
\end{theorem}

If we additionally assume that the electric potential
$V$ is semi-bounded below then we can add the bottoms of the
Neumann spectrum to the picture as was first done by
A.~Molchanov \cite{Molchanov}, though without magnetic field.
This will be important for some arguments in this paper.

\begin{theorem}\label{T:localization-N} If $V$ is semi-bounded below then
the conditions in Theorem \ref{T:localization}
are also equivalent to the following conditions:

$(f)$ $\mu(B(x,r);H_{a,V})\to +\infty$ as $x\to\infty$ for any fixed
$r>0$;

$(g)$ there exists $r>0$ such that
$\mu(B(x,r); H_{a,V})\to +\infty$ as $x\to\infty$.
\end{theorem}

Finally we can weaken the requirement on $\La$ by use of capacity
through the Molchanov condition:

\begin{theorem}\label{T:localization-M}
Let us assume that $V$ is semi-bounded below. Then the conditions
$(a)-(g)$ in Theorems  \ref{T:localization} and
\ref{T:localization-N} are equivalent
to the existence of a measurable
function $\La:\R^n\to\R$ such that

$1)$ $\La$ is semi-bounded below and satisfies $(M_c)$ with
some $c>0$;

$2)$ the operator inequality $H_{a,V}\ge\La(x)$ holds
in the same sense as in Theorem \ref{T:localization}.
\end{theorem}

For $\La=V$ and $a=0$ this gives a sufficiency part of the Molchanov
theorem. It also implies the following result from
\cite{Avron-Herbst-Simon}:

\begin{corollary}\label{C:V-to-B}
If $V$ is semi-bounded below and  $H_{0,V}$ has a discrete spectrum
$($or, equivalently, $V$ satisfies $(M_c)$ with
some $c>0$$)$,
then $H_{a,V}$ also has a discrete spectrum
for any magnetic potential~$a$.
\end{corollary}

Our proof of this statement is in fact purely variational
and it is completely different from the proof given in
\cite{Avron-Herbst-Simon} where semigroup methods are used.
However the regularity requirements for $a$ and $V$ are much
weaker in \cite{Avron-Herbst-Simon}.

\medskip
Finally let us formulate a convenient sufficient condition
which does \textit{not} require $V$ or $\La$
to be semi-bounded below.

\begin{theorem}\label{T:localization-cap}
Assume that
there exists a measurable function $\La:\R^n\to\R$ such that
the following conditions are satisfied:

$1)$ the operator inequality
\begin{equation}\label{E:Ha0-minorant}
H_{a,0}\ge\La(x)
\end{equation}
holds in the same sense as above;

$2)$ there exists $\de\in [0,1)$ such that
the effective potential
\begin{equation}\label{E:eff-pot}
\Veffde(x)=V(x)+\de\La(x)
\end{equation}
is semi-bounded below and
satisfies $(M_c)$ with some $c>0$.

Then $H_{a,V}$ is essentially self-adjoint, semi-bounded below
and has a discrete spectrum.
\end{theorem}

\subsection{Sufficient conditions ($n=2$)}\label{SS:suff2}

The case $n=2$ is much simpler than the
general case because the magnetic field does not change direction
(and may only change sign). We will identify the magnetic field
with its component $B_{12}$. Note that by changing enumeration of
coordinates we can change sign of $B=B_{12}$. Some uncertainty
principle related arguments lead to the following fact
established by J.~Avron, I.~Herbst and B.~Simon
\cite{Avron-Herbst-Simon}:

\begin{theorem}\label{T:AHS-2dim}
If $n=2$ and $|B(x)|\to\infty$ as $x\to\infty$, then
$\sigma=\sigma_d$ for $H_{a,0}$ (hence for $H_{a,V}$
with arbitrary $V\ge 0$).
\end{theorem}

\medskip
This is the simplest result which shows that magnetic field alone
may cause a localization (i.e.\ a discrete spectrum) for a quantum
particle. Classically this can be understood from the fact that
a strong magnetic field causes a fast rotation of a charged moving
particle without changing its kinetic energy, and in this way
the field impedes possible escape of the particle to infinity.

\medskip
Note that $|B(x)|\to\infty$ means in fact that either $B$ or $-B$
tend to $+\infty$ as $x\to\infty$. Hence the following theorem
improves the result above:

\begin{theorem}\label{T:2dim-0} Assume that the following conditions
are satisfied:

$(a)$ there exists $C>0$ such that $B(x)\ge -C$, $x\in \R^2$;

$(b)$ the Molchanov condition $(M_c)$ is satisfied for
$B(x)$ or, equivalently, for $|B(x)|$ $($instead of $V(x)$$)$
with some $c>0$.

Then $H_{a,0}$ has a discrete spectrum.
\end{theorem}

The simplest explicit  result which takes into account 
both electric and magnetic field is given by the following 

\begin{theorem}\label{T:simplest}
Assume that $n=2$ and there exists $\de\in [-1,1]$
such that for the effective potential
\begin{equation}\label{E:Veff}
\Veffde(x)=V(x)+\de B(x)
\end{equation}
we have $\Veffde(x)\to +\infty$ as $x\to\infty$.
Then $H_{a,V}$ is essentially self-adjoint,
semi-bounded below and has a discrete spectrum. 
\end{theorem}

V.~Ivrii noticed that the range of possible $\de$ is precise here:
the conclusion does not hold if we take any $\de\not\in [-1,1]$.

\ms
The following  theorem strengthens Theorem \ref{T:2dim-0} by
taking into account
the influence of the electric potential $V$:

\begin{theorem}\label{T:2dim} Assume that $n=2$
and there exists $\de\in (-1,1)$ such that
the effective potential $V_{eff}$ given by (\ref{E:Veff})
is semi-bounded below and  satisfies
the Molchanov condition $(M_c)$
with some $c>0$. Then $H_{a,V}$ is essentially self-adjoint,
semi-bounded below and has a discrete spectrum.
\end{theorem}

Clearly Theorem \ref{T:2dim}
implies Theorem \ref{T:2dim-0} (take $V=0$). It also implies
the sufficiency of the Molchanov condition in case $B=0$.
Note however that the conditions of the theorem do not imply any
growth of $V$ or $B$. If $|B|$ itself tends to $\infty$ as
$x\to\infty$, then $V$ is even allowed to go to $-\infty$,
though in this case $|V|$ must be dominated by
$\de |B|$ with some positive $\de<1.$

Theorem \ref{T:2dim} also implies 
Theorem \ref{T:simplest} except for the 
extreme values $\de=\pm 1$. 

\medskip
Unfortunately we are not aware whether any of the conditions
in Theorem \ref{T:2dim}
is necessary (though they are necessary if $B=0$ and $V$ is
semi-bounded below due to the Molchanov result quoted above).

Even Theorem \ref{T:simplest}, which does not require any use of capacity
and can be proved by elementary means (see Sect.\ref{SS:sufficient2}),
seems to be absent in the literature.
A stronger result without capacity can be obtained if we replace
capacity by the Lebesgue measure. Namely, let us formulate
the corresponding condition for a semi-bounded below function $V$:
there exists $r_0>0$ such that for any $r\in(0,r_0)$
\begin{equation}\label{E:Mes-cond2}
\underset{F}{\inf}
\left\{\left.\int_{B(x,r)\setminus F}V(y)dy \right|\;\mes(F)\le
c r^N\right\}
\to+\infty \ \hbox{as}\ x\to\infty.
\tag{$\tilde M_{c,N}$}
\end{equation}

It follows from a well-known estimate of measure by
capacity, that for any $c>0$ and $N>0$ the condition
$(\tilde M_{c,N})$ implies $(M_{c'})$ for some $c'>0$
(see the proof of part 2) of Theorem 6.1 in \cite{Kondrat'ev-Shubin}).
Therefore Theorem \ref{T:2dim} implies the following

\begin{corollary}\label{C:2dim-mes}
Assume that there exists $\de\in(-1,1)$ such that
the effective potential $\Veffde$ given by
(\ref{E:Veff}) is semi-bounded below and satisfies $(\tilde M_{c,N})$
with some $c>0$ and $N>0$.
Then $H_{a,V}$ is essentially self-adjoint, semi-bounded below
and has a discrete spectrum.
\end{corollary}

Finally note that an elementary argument given in
the proof of Corollary 6.2 in \cite{Kondrat'ev-Shubin}
gives a sufficient condition which is stronger than
(\ref{E:Mes-cond2}) but very easy to check.

\begin{corollary}\label{C:2dim-elem}
Assume that there exists $\de\in (-1,1)$ such that
$V_{eff}^{(\de)}$ is semi-bounded below and
for any $A>0$ and any small $r>0$
\begin{equation}\label{E:2dim-elem1}
\mes\{y|\;y\in B(x,r),\ V_{eff}^{(\de)}(y) \le A\}\to 0
\quad \hbox{as} \quad x\to\infty.
\end{equation}
Then $H_{a,V}$ is essentially self-adjoint, semi-bounded below
and has a discrete spectrum.
\end{corollary}

\subsection{Sufficient conditions ($n\ge 3$)}\label{SS:suff3}

The behavior of the spectrum of the magnetic Schr\"odinger
operator in dimensions $n\ge 3$ is much more complicated than in
dimension 2 because of possible varying direction of $B$.
In particular none of the results formulated above for $n=2$
holds for $n\ge 3$. A.~Dufresnoy \cite{Dufresnoy} gave the first
example of the operator $H_{a,0}$ with
\begin{equation}\label{E:B-to-infty}
|B(x)|\to\infty \quad \hbox{as} \quad x\to\infty,
\end{equation}
and yet with non-compact resolvent (or,
equivalently, with non-discrete spectrum) in an arbitrary
dimension $n\ge 3$. J.~Avron, I.~Herbst and B.~Simon
\cite{Avron-Herbst-Simon} gave
sufficient conditions for
the discreteness of the spectrum of $H_{a,0}$ which in addition to
(\ref{E:B-to-infty}) require that the direction of $B$ varies
sufficiently slowly. A more explicit condition of this kind
(in terms of estimates for derivatives of the direction of $B$)
was given in \cite{Dufresnoy}. Later A.~Iwatsuka \cite{Iwatsuka1}
improved these results, giving almost precise estimates of this kind.
He also produced a series of spectacular examples. One of them
shows that no growth condition for $|B(x)|$ (i.e.\ a condition of the
form $|B(x)|\ge \rho(x)$ with a fixed continuous function $\rho$) would
suffice for the discreteness of the spectrum of $H_{a,0}$. This
is in drastic contrast with the results for $n=2$ formulated
above. Another example given in \cite{Iwatsuka1} shows that
the condition (\ref{E:B-to-infty}) is also not necessary for the
discreteness of spectrum of $H_{a,0}$. (This is of course
less surprising because you can expect that integrally small
perturbations of $B$ should not generally affect the discreteness of
spectrum.)

\medskip
However we will show in Sect.\ref{S:sufficient} that 
some explicit sufficient conditions
for the discreteness of spectrum can still be formulated in terms
of effective potentials similarly to the case $n=2$. 
The appropriate effective potentials will include both electric and
magnetic fields. They will incorporate information
about the direction of the magnetic field. The above mentioned results
of  J.~Avron, I.~Herbst and B.~Simon,
A.~Dufresnoy  and A.~Iwatsuka will follow if we impose additional conditions
on the direction of the magnetic field. These conditions allow us to simplify 
the form of the effective potential.

\ms
 Here we will formulate just one result coming from this approach 
(others can be found in Sect.\ref{S:sufficient}). 

\ms By $\Lip(X)$ we will denote the set of all Lipschitz functions
on any metric space. We will only use $X$ which are locally compact.
In this case $\Liploc(X)$ will denote the space of functions which are
locally Lipschitz on $X$. Recall that Lipschitz functions on an open subset 
in $\R^n$ are exactly the ones which have bounded distributional derivatives
(see e.g.\ \cite{Mazya8}).

\ms
Let us define a smoothed direction of the magnetic field as follows:
\begin{equation}\label{E:Duf-mod}
A_{jk}(x)=\chi(|B(x)|)\frac{B_{jk}(x)}{|B(x)|},
\end{equation}
where $\chi\in \Lip([0,\infty))$, $\chi(r)=0$ if $r\le 1/2$,
$\chi(r)=1$ if $r\ge 1$, $\chi(r)=2r-1$ if $1/2\le r\le 1$, so
$0\le \chi(r)\le 1$ and $|\chi'(r)|\le 2$
for all $r$.

\begin{theorem}\label{T:Duf-gen-maj}
Let us assume that $B_{jk}\in \Lip(\R^n)$ for all $j,k$;
$A_{jk}$ are defined by (\ref{E:Duf-mod}),
and a positive measurable function
$X(x)$ in $\R^n$ satisfies
\begin{equation}\label{E:A-majorant}
\sum_{k=1}^n\left|\frac{\pa A_{kj}(x)}{\pa x^k}\right|\le
X(x),\quad x\in\R^n,
\end{equation}
for all $j$. Then for any $\eps>0$ and $\de\in [0,1)$ define an effective
potential
\begin{equation}\label{E:Veff-de-gen}
\Veffdeeps(x)=V(x)+\frac{\de}{n-1+\eps}|B(x)|-
\frac{n\de}{4\eps(n-1+\eps)}X^2(x).
\end{equation}
If there exist $\eps>0$ and $\de\in [0,1)$, such that $\Veffdeeps$
is bounded below and satisfies the Molchanov condition $(M_c)$
with some $c>0$ (in particular, this holds if
$\Veffdeeps(x)\to +\infty$ as $x\to\infty$),
then $H_{a,V}$ is essentially self-adjoint,
semi-bounded below and has a discrete spectrum.
\end{theorem}

Note that in case $B=0$ and $V$ semi-bounded below the condition on 
$V$ in this theorem becomes necessary and sufficient due 
to the Molchanov theorem.

\medskip
In Section \ref{S:sufficient} we will establish that
imposing different regularity conditions on $B$, e.g.\ as in
\cite{Avron-Herbst-Simon, Dufresnoy, Iwatsuka1},
leads to different types of effective potentials so that
results similar to Theorem \ref{T:2dim} hold.

Here we will give the simplest example of this kind.
As explained above it is natural to impose some a priori
conditions of regularity on $B$.
We will mainly formulate them in the form
of estimates for derivatives of $B$.
We will always assume
that $B\in \Liploc$, i.e.\ $B_{jk}\in \Liploc$ for all $j,k$.
Following A.~Iwatsuka \cite{Iwatsuka1}
we will use the estimates of the form

\begin{equation}\label{E:B-der-alpha}
|\nabla B(x)|\le C(1+|B(x)|)^\al,\quad x\in\R^n,
\tag{$B_\al$}
\end{equation}
where $\al>0$, $C>0$ and $\nabla B$ means vector whose components
are all possible first order derivatives $\pa B_{jk}/\pa x^l$.
We will write that $B$ satisfies (\ref{E:B-der-alpha}) if there
exists $C>0$ such that this estimate is satisfied with the given
$\al$. A little bit stronger condition used in \cite{Dufresnoy}
and \cite{Iwatsuka1} is
\begin{equation}\label{E:B-der-alpha0}
|\nabla B(x)|(1+|B(x)|)^{-\al}\to 0 \quad\hbox{as}\quad
x\to\infty.
\tag{$B_\al^0$}
\end{equation}

A.~Dufresnoy \cite{Dufresnoy} proved that the conditions
(\ref{E:B-to-infty}) and $(B_{3/2}^0)$ imply that the spectrum of
$H_{a,0}$ is discrete. In fact he proved that instead of $(B_{3/2}^0)$
it is sufficient to require a slightly weaker condition
\begin{equation}\label{E:B-der-dir}
\left|\nabla\left(\frac{B(x)}{|B(x)|}\right)\right|\cdot
|B(x)|^{-1/2}\to 0 \quad\hbox{as}\quad x\to\infty.
\notag
\end{equation}
A.~Iwatsuka \cite{Iwatsuka1} proved that in fact $(B_2^0)$ together
with (\ref{E:B-to-infty}) imply the discreteness of spectrum for
$H_{a,0}$. He also provided an example which shows that $(B_2^0)$ can not be
replaced by $(B_2)$ (hence it can not be replaced by $(B_\al)$ or $(B_\al^0)$
with any $\al>2$).

The following theorem compared with the above mentioned results
takes into account the behavior of $V$.

\begin{theorem}\label{T:B-and-V}
Assume that

$(a)$ $B$ satisfies $(B_{3/2}^0)$,

\medskip
and

\medskip
$(b)$ there exists $\de\in [0,1)$  such that
the effective potential
\begin{equation}\label{E:Veff3}
\Veffde=V+\frac{\de}{n-1}|B|
\end{equation}
is semi-bounded below
and satisfies the Molchanov condition $(M_c)$ with a small $c>0$
(possibly depending on $\de$).

\medskip
Then $H_{a,V}$ is essentially self-adjoint, semi-bounded below and
has a discrete spectrum.
\end{theorem}

\begin{corollary}\label{C:3dim}
Assume that

\medskip
$(a)$ $B$ satisfies $(B_{3/2}^0)$

\medskip
and

\medskip
$(b)$ there exists $\de\in[0,1)$ such that
$\Veffde(x)\to +\infty$ as $x\to\infty$
for the effective potential $\Veffde$
defined by (\ref{E:Veff3}).

\medskip
Then $H_{a,V}$ is essentially self-adjoint, semi-bounded below and
has a discrete spectrum.
\end{corollary}

Due to the arguments given in Sect.6.1 of \cite{Kondrat'ev-Shubin}
it is also possible to replace capacity by the Lebesgue measure.
Namely, let us introduce the corresponding condition for  a
semi-bounded below function $V$:

\begin{equation}\label{E:Mes-cond3}
\underset{F}{\inf}
\left\{\left.\int_{B(x,r)\setminus F}V(y)dy \right|\;\mes(F)\le
c r^n\right\}
\to+\infty \ \hbox{as}\ x\to\infty.
\tag{$\tilde M_{c}$}
\end{equation}

As shown in Sect.6.1 of \cite{Kondrat'ev-Shubin}, for any $c>0$
there exists $c'>0$  such that $(\tilde M_c)$ implies $(M_{c'})$.
Therefore we have the following

\begin{corollary}\label{C:3dim-mes}
Assume that $B$ satisfies $(B_{3/2}^0)$ and

\medskip
$(b)$ there exists $\de\in [0,1)$ such that the
the effective potential (\ref{E:Veff3})
is semi-bounded below and
 there exists $c>0$  such that
the condition (\ref{E:Mes-cond3}) holds for $\Veffde$
(instead of $V$);

\medskip
Then $H_{a,V}$ is essentially self-adjoint, semi-bounded below and
has a discrete spectrum.
\end{corollary}

Finally let us formulate an analogue of Corollary \ref{C:2dim-elem}
for $n\ge 3$.

\begin{corollary}\label{C:3dim-elem}
Assume that $B$ satisfies $(B_{3/2}^0)$ and

\ms
$(b)$ there exists $\de\in[0,1)$ such that
the effective potential (\ref{E:Veff3})
is semi-bounded below and for any $A>0$
and any small $r>0$
\begin{equation}\label{E:3dim-elem2}
\mes\{y|\;y\in B(x,r),\ \Veffde(y)\le A\}\to 0
\quad \hbox{as} \quad x\to\infty.
\end{equation}

\medskip
Then $H_{a,V}$ is essentially self-adjoint, semi-bounded below and
has a discrete spectrum.
\end{corollary}

Note that if we assume that $V$ is semi-bounded below,
then we can always take $V+|B|$ as an effective potential in all
statements above.

\subsection{Acknowledgments}

Most part of this work was done during our stay in the program
Research in Pairs of Oberwolfach Forschungsinstitut f\"ur
Mathematik in May-June 1999. We gratefully acknowledge the
generous support of this Institute and its kind and helpful
personnel.

The second author was  also partially supported by NSF
grant DMS-9706038.

The authors are grateful to M.S.~Birman, J.-M.~Bismut, V.~Ivrii, J.~Lott,  
V.G.~Maz'ya and M.Z.~Solomyak for useful discussions. 

\section{Essential self-adjointness}\label{S:self}

The goal of this section is to give a simple proof of the following
simplest version of the  theorem on essential
self-adjointness of any semi-bounded magnetic Schr\"odinger
operator (see \cite{Povzner, Wienholtz, Glazman, Simader78, Shubin00} 
for other versions). 

\begin{theorem}\label{T:self}
 Assume that $a_j\in \Liploc(\R^n)$,
$V\in L^\infty_{loc}(\R^n)$ and the corresponding magnetic
Schr\"odinger operator $H_{a,V}$ is semi-bounded below on $C_c^\infty(\R^n)$,
i.e.\ there exists a constant $C\in\R$ such that
\begin{equation}\label{E:Hsemibound}
(H_{a,V}u,u)\ge -C(u,u),\quad u\in C_c^\infty(\R^n).
\end{equation}
Then $H_{a,V}$ is essentially self-adjoint.
\end{theorem}

\textbf{Proof.} We will extend the Wienholtz proof of the
Povzner theorem as it is explained by I.M.~Glazman \cite{Glazman}.

Let us recall that if $g\in\Lip(\Om)$ where $\Om$ is an open
subset in $\R^n$, then $\pa g/\pa x^j\in L^\infty(\Om)$,
$j=1,\dots,n$  (see e.g.\ \cite{Mazya8}, Sect.1.1), where the
derivatives are understood in the sense of distributions
(but also exist almost everywhere). This implies that the operator
$H_{a,V}$ is well defined on $C_c^\infty(\R^n)$
(and maps this space into $L^2(\R^n)$)
as well as on $L^2(\R^n)$ (which it maps to the space
of distributions on $\R^n$).

Note that adding $(C+1)I$ to $H_{a,V}$
we can assume that $H_{a,V}\ge I$ on $C_c^\infty(\R^n)$,
i.e.
\begin{equation}\label{E:semi-I}
h_{a,V}(u,u)\ge (u,u),\quad u\in C_c^\infty(\R^n).
\end{equation}
If this is true then it is well known (see e.g.\ \cite{Glazman})
that the essential self-adjointness of $H_{a,V}$ is equivalent
to the fact that the equation
\begin{equation}\label{E:ker}
H_{a,V}u=0
\end{equation}
has no non-trivial solutions in $L^2(\R^n)$ (understood in the sense of
distributions).

Assume that $u$ is such a solution. First note that it is in
$W^{2,2}_{loc}(\R^n)$ (i.e.\ has distributional derivatives of order
$\le 2$ in $L^2_{loc}(\R^n)$) due to a simple elliptic regularity argument
(see e.g.\ Lemma 4.1 in \cite{Shubin3} for more details).

Let us take a cut-off function $\phi_R\in C_c^\infty(\R^n)$
with the following properties:
\begin{eqnarray}
&0\le\phi_R\le 1;
\notag\\
&\phi_R=1 \quad \textit{on}\quad B(0,R) \quad
\textit{and} \quad 0\quad \textit{on}\quad
\R^n\setminus B(0,2R);
\notag\\
&\eps_R:=\underset{x\in\R^n}{\sup}|\nabla\phi_R(x)|\to 0
\quad \textit{as} \quad R\to\infty.
\notag
\end{eqnarray}

Then denoting $u_R=\phi_R u$ we see that $u_R$ is in the domain
of the minimal operator associated with $H_{a,V}$, hence
\begin{equation}\label{E:semi-cut}
\|u_R\|^2\le (H_{a,V}u_R,u_R).
\end{equation}

Let us calculate $H_{a,V}u_R$ using the Leibniz
type formula for $P_j$:
\begin{equation}\label{E:Leibniz}
P_j(fg)=(P_jf)g+f(D_jg),
\end{equation}
where $D_j=-i\pa/\pa x^j$. Applying this formula twice
to calculate $P_j^2(\phi_Ru)$ and summing up we
easily obtain due to (\ref{E:ker}):
\begin{eqnarray}\label{E:HaV-cut}
&H_{a,V}u_R&=\phi_R H_{a,V}u-2(\nabla\phi_R)\cdot(\nabla_au)
-u\De\phi_R\\
&&=-2(\nabla\phi_R)\cdot(\nabla_au)
-u\De\phi_R \notag\\
&&=-2\nabla\phi_R\cdot\nabla u-2i(a\cdot\nabla\phi_R)-u\De\phi_R.
\notag
\end{eqnarray}
Therefore due to \ref{E:semi-cut} we have
\begin{equation}
\|\phi_Ru\|^2\le \int_{\R^n}
(-2\nabla\phi_R\cdot\nabla u-2i(a\cdot\nabla\phi_R)u-u\De\phi_R)
\phi_R\bar u dx.
\end{equation}
Since $(H_{a,V}u_R,u_R)$ is real, we can replace
the right hand side here by the complex conjugate expression.
Adding the two estimates obtained in such a way
and dividing by 2 we see
that the term with the magnetic potential $a$ cancels and we get,
applying integration by parts:
\begin{eqnarray}
&\|\phi_R u\|^2&\le \int_{\R^n}\left[-\phi_R(\nabla\phi_R)\cdot
(\bar u\nabla u+u\nabla\bar u)-\phi_R(\De\phi_R)|u|^2\right]dx
\notag \\
&&=\int_{\R^n}\left[-\phi_R(\nabla\phi_R)\cdot
\nabla(|u|^2)-\phi_R(\De\phi_R)|u|^2\right]dx
\notag\\
&&=\int_{\R^n}\left[\phi_R(\De\phi_R)|u|^2
+|\nabla\phi_R|^2|u|^2-\phi_R(\De\phi_R)|u|^2\right]dx
\notag\\
&&=\int_{\R^n}|\nabla\phi_R|^2|u|^2 dx.
\notag
\end{eqnarray}
In particular we obtain using the conditions on $\phi_R$ above:
\begin{equation}
\int_{B(0,R)}|u|^2 dx\le \eps_R^2\int_{B(0,2R)}|u|^2dx.
\notag
\end{equation}
Allowing $R$ to go to $+\infty$ we see that $\|u\|^2=0$, hence
$u\equiv 0$.
$\ecarre$

\medskip
\textbf{Remark.} The local condition $V\in L^\infty_{loc}(\R^n)$ can be 
considerably weakened. For example, it is sufficient to require that $V=V_++V_-$,
where $V_+\ge 0$, $V_+\in L^2_{loc}(\R^n)$, $V_-\le 0$, $V_-\in L^p_{loc}(\R^n)$
with $p=2$ if $n\le 3$, $p>2$ if $n=4$, and $p=n/2$ if $n\ge 5$. 
(See e.g.\ \cite{Simader78, Shubin00}.)

\section{Localization}\label{S:localization}

In this section we will prove different localization theorems
which were formulated in Sect.\ref{SS:localization}
and provide an important preliminary material
related to compactness arguments and estimates of the bottoms of
the Dirichlet and Neumann spectra.

We will use notations from previous sections.
In particular, $a,V$ will always denote magnetic and
electric potential with the same regularity as in Section
\ref{S:Intro}.

We will start with the following elementary and well known
diamagnetic inequality (see e.g.\ \cite{Kato2, Simon76, Lieb-Loss}):

\begin{lemma}\label{L:grad-mod-est}
Let $a$ be an arbitrary magnetic potential
(with components from $C^1$).
Let $u$ be a complex valued Lipschitz function in an open
set $U\subset\R^n$. Then $|u|$ is also Lipschitz and
\begin{equation}\label{E:grad-mod-est}
|\nabla |u||\le |\nabla_a u| \quad \textit{a.e.,}
\end{equation}
where a.e.\ means almost everywhere with respect to
the Lebesgue measure.
\end{lemma}

\medskip
Let us assume that $H_{a,V}$ is bounded below, hence
essentially self-adjoint due to Theorem \ref{T:self}.
Without loss of generality we can assume hereafter that $H_{a,V}\ge I$
(or that the estimate (\ref{E:semi-I}) is satisfied).

The essential self-adjointness of $H_{a,V}$ means that
$C_c^\infty(\R^n)$ is its \textit{core}, i.e.\ the closure of
$H_{a,V}$ from the initial domain $C_c^\infty(\R^n)$
is a self-adjoint operator in $L^2(\R^n)$. It follows
that $C_c^\infty(\R^n)$ is also a core for the corresponding
quadratic form.

Denote
\begin{equation}\label{E:L}
\cL =\{u\in C_c^\infty(\R^n)|\;h_{a,V}(u,u)\le 1\}.
\end{equation}

\begin{lemma}\label{L:compactness}
$\sigma=\sigma_d$ if and only if $\cL$ is precompact in $L^2(\R^n)$.
\end{lemma}

\textbf{Proof.} The proof is the same as the proof of
Lemma 2.2 in \cite{Kondrat'ev-Shubin} and it is essentially
abstract. Clearly $\sigma=\sigma_d$ is equivalent to saying
that for $H=H_{a,V}$ we have
\begin{equation}\label{E:Dom-quad}
\{u|\;u\in \Dom(H^{1/2}),\ \|H^{1/2}u\|\le 1\}\
\end{equation}
is precompact in $L^2(\R^n)$. Here $\Dom(H^{1/2})$ also coincide
with the domain of the quadratic form which is the closure of
$h_{a,V}$. Since $C_c^\infty(\R^n)$ is a core for the quadratic
form too, we see that precompactness of the set (\ref{E:Dom-quad})
is equivalent to the precompactness of $\cL$. $\ecarre$

\begin{lemma}\label{L:smallintegral} {\sc({Small Tails Lemma})}
Let  us assume as above that
$H_{a,V}$ is essentially self-adjoint and semi-bounded below
so that (\ref{E:semi-I}) holds. Then
$\sigma=\sigma_d$ if and only if the following small tails
condition is satisfied: for any $\eps>0$
there exists $R>0$ such that
\begin{equation}\label{E:M-B}
\int_{\R^n\setminus B(0,R)}|u|^2dx<\eps\ \mbox{for any}\ u\in \cL,
\end{equation}
or, in other words,
$$
\int_{\R^n\setminus B(0,R)}|u|^2dx\to 0\ \mbox{as}\ R\to\infty,\
\mbox{uniformly in} \ u\in \cL\;.
$$
\end{lemma}

\textbf{Proof.}
Again the proof is similar to the proof of Lemma
2.3 in \cite{Kondrat'ev-Shubin}, though a small additional
argument is needed to avoid using semi-boundedness of $V$
which was used in \cite{Kondrat'ev-Shubin}.

Clearly $\sigma=\sigma_d$ (or precompactness of $\cL$) implies
the small tails condition because any precompact set has a
$\eps$-net for any $\eps>0$.

Vice versa, assume that the small tails condition is fulfilled.
Then the precompactness of $\cL$ would be equivalent to the
precompactness of any restriction
$$
\cL_R=\left\{\left. u|_{B(0,R)}\right|\; u\in\cL\right\},
$$
where $R>0$. Note that the condition (\ref{E:semi-I}) implies that
the set $\cL_R$ is bounded in $L^2(B(0,R))$. But then by the
Sobolev--Kondrashov compactness theorem it is sufficient
to establish a uniform $L^2(B(0,R))$-boundedness of the gradients of
functions $u\in\cL_R$ (for any fixed $R>0$). Since we assume
that $a_j\in L^\infty_{loc}$, it is sufficient to establish that
the magnetic gradients $\nabla_au$ are uniformly bounded in
$L^2(B(0,R))$.  This in turn follows from the definition of
$\cL$ and uniform $L^2$-boundedness of the functions $u\in \cL$
combined with the local boundedness of the potential $V$.
$\ecarre$

\medskip
\textbf{Remark.} The requirement $u\in C_c^\infty(\R^n)$ in
the definition of $\cL$ (and in
Lemmas \ref{L:compactness} and \ref{L:smallintegral}) can be
replaced by the requirement $u\in \Lip_c(\R^n)$
(the set of all Lipschitz functions with
compact support in $\R^n$)
because the space $\Lip_c(\R^n)$ is intermediate between
$C_c^\infty(\R^n)$ and $\Dom\left(H_{a,V}^{1/2}\right)$.

\medskip
Now let us take a covering of $\R^n$ by balls $B(x_k,r)$, $k=1,2,\dots$,
of fixed radius $r>0$, such that this covering has a finite
multiplicity. Now take a  partition of unity on $\R^n$ consisting of functions
$e_k\in C^\infty_c(\R^n)$, $k=1,2,\dots$,  such that
$0\le e_k\le 1$, $\supp e_k\subset B(x_k,r)$, and
\begin{equation}\label{E:part-unity}
\sum_{k=1}^\infty e_k^2=1,
\notag
\end{equation}
\begin{equation}\label{E:der-part-unity}
|\nabla e_k|\le C, \quad k=1,2,\dots,
\notag
\end{equation}
where $C$ does not depend on $k$.

The main tool in proving the localization theorem is the following
IMS localization formula
\begin{equation}\label{E:IMS-localization}
H_{a,V}=\sum_{k=1}^\infty J_kH_{a,V}J_k -
\sum_{k=1}^\infty |\nabla e_k|^2,
\end{equation}
where $J_k$ is the multiplication operator by $e_k$ in
$L^2(\R^n)$.
Proofs of different versions of this formula can be
found in \cite{Simon83}, \cite{Cycon-Froese-Kirsch-Simon}
(Sect. 3.1), \cite{Shubin96}.  Formally only the  case $a=0$
is treated in \cite{Simon83, Cycon-Froese-Kirsch-Simon},
though the proof works with arbitrary $a$. Much more general case
(second order differential operators on manifolds) is considered in
\cite{Shubin96}, Sect.3.

\medskip
Now we will fix an operator $H_{a,V}$ and denote for brevity
$\la(B(x,r))=\la(B(x,r);H_{a,V})$,
$\mu(B(x,r))=\mu(B(x,r);H_{a,V})$.

\medskip
\textbf{Proof of Theorem \ref{T:localization}.}
$(a)\Longrightarrow (b).$ Let us assume that $(a)$ is satisfied and
fix an arbitrary $r>0$.
According to Lemma \ref{L:smallintegral}
for any $\eps>0$ there exists $R>0$ such that
$|x|>R$ implies that
$$
\int_{B(x,r)} |u(x)|^2 dx<\eps,
$$
as soon as $u\in C_c^\infty(B(x,r))$ and $h_{a,V}(u,u)\le 1$.
It follows that $\la(B(x,r))\ge 1/\eps$, which implies $(b)$.

\medskip
Clearly $(b)\Longrightarrow (c)$.

\medskip
$(c)\Longrightarrow (d)$. Let us fix $r>0$ such that
$(c)$ is satisfied and choose a covering of $\R^n$
by the balls $B(x_k,r)$, $k=1,2,\dots,$, and a partition
of unity $\{e_k|\;k=1,2,\dots,\}$ with $\supp e_k\subset B(x_k,r)$
and with the properties formulated above. Then for any
$u\in C_c^\infty(\R^n)$ we obtain from (\ref{E:IMS-localization}):
\begin{equation}
h_{a,V}(u,u)=\sum_k h_{a,V}(e_ku, e_ku)-\sum_k|\nabla e_k|^2|u|^2
\ge (\La u,u),
\notag
\end{equation}
where
\begin{equation}
\La(x)=\sum_k \la(B(x_k,r))e_k^2(x)-\sum_k |\nabla e_k|^2.
\notag
\end{equation}
The first sum tends to $+\infty$ as $x\to\infty$ due to the
condition $(c)$, whereas the second sum is bounded. This implies
that $\La(x)\to +\infty$ as $x\to\infty$, hence $(d)$ is
fulfilled.

\medskip
The implication $(d)\Longrightarrow (e)$ is obvious.

\medskip
$(e)\Longrightarrow (a).$
Assume that $(e)$ is satisfied with the function $\La(x)$. Let us
take $u\in\cal{L}$. Then
\begin{equation}
\int_{\R^n}\La(x)|u(x)|^2dx\le 1.
\notag
\end{equation}
Therefore,
\begin{equation}
\int_{|x|\ge R}|u(x)|^2dx \le
\left(\underset{\{x|\;|x|\ge R\}}{\inf} \{\La(x)\}\right)^{-1}
\to 0\quad\hbox{as}\
R\to\infty,
\notag
\end{equation}
so $(a)$ follows from Lemma \ref{L:smallintegral}.
$\ecarre$

\medskip
Now we assume that $V\ge 0$ and proceed to some preparatory
material which is needed to prove Theorem \ref{T:localization-N}.

\medskip
Denote temporarily $B_r=B(0,r)\subset\R^n$ and define
\begin{equation}\notag
\|\psi\|=\|\psi\|_{L^2(B_r)}=
\left(\int_{B_r}|\psi|^2dx\right)^{1/2},
\quad \|\psi\|_t=\|\psi\|_{L^2(B_{tr})},
\end{equation}
where $0<t\le 1$.
Similarly define

\begin{eqnarray*}
&&\|\nabla\psi\|=\|\nabla\psi\|_{L^2(B_r)}=
\left(\int_{B_r}|\nabla\psi|^2dx\right)^{1/2},\\
&&\|\nabla\psi\|_t=\|\nabla\psi\|_{L^2(B_{tr})}.
\end{eqnarray*}

\begin{lemma}\label{L:tball}
The following estimates hold true for any magnetic potential $a$:

\begin{equation}\label{E:smallerball1}
\|\psi\|\le t^{-n/2}\|\psi\|_t+2^{n+1}r(1-t)\|\nabla_a\psi\|, \quad
t\in [1/2,1];
\end{equation}

\begin{equation}\label{E:smallerball2}
\|\psi\|^2\le 2t^{-n}\|\psi\|^2_t+2^{2n+3}r^2(1-t)^2\|\nabla_a\psi\|^2,
\quad t\in [1/2,1].
\end{equation}
In these estimates we assume that $\psi\in \Lip(B_r)$.
\end{lemma}

\textbf{Proof.} With $a=0$ this Lemma was proved in
\cite{Molchanov} (for cubes instead of balls) and in
\cite{Kondrat'ev-Shubin}(see Lemma 2.8 there). Applying this
particular case to $|\psi|$ and using Lemma \ref{L:grad-mod-est},
we obtain the desired result.
$\ecarre$

\medskip
The following Lemma for the case $a=0$ was proved in
\cite{Molchanov} (for cubes) and \cite{Kondrat'ev-Shubin} (see Lemma 2.9
there):

\begin{lemma}\label{L:comparison} Let us assume that $V\ge 0$.
Then
\begin{equation}\label{E:mu-lambda-mu}
\mu(B(x,r))\le \la(B(x,r))\le C_1 \mu(B(x,r))+ C_2 r^{-2},
\end{equation}
for any $x\in \R^n$, and any $r>0$. Here $C_1$ and $C_2$
depend only on $n$; for example we can take
$C_1=2^{n+3}(1+2^{2n+6})$ and $C_2=2^{n+7}$.
\end{lemma}

\textbf{Proof.} The proof given in \cite{Kondrat'ev-Shubin} works
in our case (with arbitrary $a$) if we replace $\nabla$ by
$\nabla_a$,
use the Leibniz rule (\ref{E:Leibniz})
and apply Lemma \ref{L:tball} above. $\ecarre$

\medskip
\textbf{Proof of Theorem \ref{T:localization-N}.}
Equivalence of $(b)$ (respectively $(c)$) from Theorem \ref{T:localization}
to $(f)$ (resp. $(g)$) from Theorem \ref{T:localization-N}
follows from Lemma \ref{L:comparison} provided we additionally
assume that $V\ge 0$. $\ecarre$

\medskip
\textbf{Proof of Theorem \ref{T:localization-M}.}
Without loss of generality we can assume that $V\ge 0$.
The necessity part of the theorem follows from
Theorem \ref{T:localization}.
(If $\sigma=\sigma_d$ then we can even make $\La(x)\to +\infty$
as $x\to\infty$.)

Now let us assume that there exists $\La(x)$ which is
semi-bounded below, satisfies $(M_{c})$ with some $c>0$,
and $H_{a,V}\ge \La(x)$, i.e.\
$$
h_{a,V}(u,u)\ge (\La|u|,|u|)
$$
for any $u\in\Lip_c(\R^n)$.
Using Lemma \ref{L:grad-mod-est} we also obtain
$$
h_{a,V}(u,u)\ge h_{a,0}(u,u)\ge h_{0,0}(|u|,|u|).
$$
Adding these two inequalities we get
$$
2h_{a,V}(u,u)\ge h_{0,\La}(|u|,|u|).
$$
Let us introduce the set
$$
\tilde\cL=\{u|\,u\in\Lip_c(\R^n),\ h_{0,\La}(u,u)\le 1\},
$$
which has the small tails property (\ref{E:M-B})
due to the Molchanov theorem and the
Small Tails Lemma \ref{L:smallintegral}.
With $\cL$ defined by (\ref{E:L}) we see that the map $u\mapsto |u|$
maps $\cL$ into $\sqrt{2}\tilde\cL$. Hence $\cL$ also has the small tails
property and we conclude from Lemma \ref{L:smallintegral} that $H_{a,V}$
has a discrete spectrum.
$\ecarre$

\medskip
\textbf{Proof of Theorem \ref{T:localization-cap}.}
We immediately conclude from the conditions that
for any function $u\in\Lip_c(\R^n)$
\begin{eqnarray*}
&h_{a,V}(u,u)&=h_{a,0}(u,u)+(Vu,u)\\
&&=(1-\de)h_{a,0}(u,u)+\de h_{a,0}(u,u)+(V|u|,|u|)\\
&&\ge (1-\de)h_{0,0}(|u|,|u|)+((\de\La+V)|u|,|u|)\\
&&=(1-\de)\left[h_{0,0}(|u|,|u|)+
\left((1-\de)^{-1}(V+\de\La)|u|,|u|\right)\right].
\end{eqnarray*}
 Now the same argument based on the Small Tails Lemma
 \ref{L:smallintegral}, as in the proof of Theorem
 \ref{T:localization-M},
 ends the proof of Theorem \ref{T:localization-cap}.
 $\ecarre$

\section{Bounded magnetic field perturbations}\label{S:bound-magn}

The main goal of this section is the proof of the following

\begin{theorem}\label{T:bound-magn}
Suppose we are given an electric potential
$V\in L^\infty_{\textit{loc}}(\R^n)$ which is
semi-bounded below, and two magnetic potentials
$a,\ta\in\Liploc(\R^n)$
such that for the magnetic field $\tB$, associated with $\ta$,
we have $\tB\in\Liploc(\R^n)$ and
\begin{equation}\label{E:B-bounded}
|\tB(x)|\le C, \quad x\in\R^n.
\end{equation}
Then  $H_{a+\ta,V}$ is essentially self-adjoint, semi-bounded below 
and has a discrete  spectrum 
if and only if this is true  for $H_{a,V}$.
\end{theorem}


\ms
We will start with the following version of Poincar\'e Lemma
which is similar to the one used by A.~Iwatsuka
\cite{Iwatsuka1}, Proposition 3.2.

\begin{lemma}\label{L:Poincare}
Let $B=\sum_{j<k}B_{jk}dx^j\wedge dx^k$ be a closed 2-form
in $B(x_0,r)\subset\R^n$ with
$B_{jk}\in \Lip(B(x_0,r))$ and $\|B_{jk}\|_\infty\le C$.
Then there exists a 1-form $a=\sum_j a_jdx^j$ with $da=B$, such that
$a_j\in\Lip(B(x_0,r))$
and $\|a_j\|_\infty\le nC$. Here $\|\cdot\|_\infty$ means the
$L^\infty$-norm on $B(x_0,r)$.
\end{lemma}

\textbf{Proof.} We can obviously assume that $x_0=0$.  Then
we can produce $a_j$ by the following explicit
formulas (see e.g.\ \cite{Warner}, p. 155--156):
$$
a_j(x)=\sum_{k=1}^n x_k\int_0^1 tB_{kj}(tx)dt,
$$
and all necessary estimates obviously follow. $\ecarre$

\ms
\textbf{Proof of Theorem \ref{T:bound-magn}.} Let us assume that
$H_{a+\ta,V}$ is essentially self-adjoint, semi-bounded below and has a discrete spectrum.
We will prove that the same holds for holds for $H_{a,V}$. 
Clearly this is sufficient to prove the theorem.

\ms
We can also assume that $V\ge 0$.

\ms
Let us choose an arbitrary ball $B(x_0,r)$. We would like
to estimate the bottom of the Dirichlet spectrum of $H_{a+\ta,V}$
which we denote, as before, by $\la(B(x_0,r);H_{a+\ta,V})$.
Using the gauge invariance we can then arbitrarily change $\ta$
in the ball $B(x_0,r)$ as soon as the relation $d\ta=\tB$
is preserved. Therefore we can use Lemma \ref{L:Poincare}
and assume that $\|\ta\|_\infty\le nrC$, where $C$ is
 the constant in \eqref{E:B-bounded} and
$$
\|\ta\|_\infty=\max_j \|\ta_j\|_\infty.
$$
Denote
$$
P'_j=\frac{1}{i}\frac{\pa}{\pa x^j}+a_j+\ta_j,\quad
P_j=\frac{1}{i}\frac{\pa}{\pa x^j}+a_j.
$$
Then for any $u\in\Lip_c(B(x_0,r))$ we have
$$
\|P'_ju\|^2=\|P_ju+\ta_j u\|^2=\|P_ju\|^2+\|\ta_ju\|^2+
2{\textrm Re} (P_ju,\ta_ju),
$$
hence
$$
h_{a+\ta,0}(u,u)=\sum_{j=1}^n \|P'_ju\|^2\le
h_{a,0}(u,u)+n\|\ta\|_\infty^2\|u\|^2+
2\sum_{j=1}^n {\textrm Re} (P_ju,\ta_ju).
$$
We have for any $\eps>0$
\begin{align*}
&2{\textrm Re} (P_ju,\ta_ju)\le 2\|P_ju\|\|\ta_ju\|\\
&\le \eps\|P_ju\|^2+\frac{1}{\eps}\|\ta_ju\|^2
\le \eps\|P_ju\|^2+\frac{1}{\eps}\|\ta\|^2_\infty\|u\|^2.
\end{align*}
Combining this with the previous estimate we obtain:
$$
h_{a+\ta,0}(u,u)\le (1+\eps)h_{a,0}(u,u)+
n\left(1+\frac{1}{\eps}\right)\|\ta\|^2_\infty(u,u).
$$
Now adding $(Vu,u)$ to the left hand side of this inequality
and $(1+\eps)(Vu,u)$ to the right hand side, we obtain
the following operator inequality
(which holds in the sense of quadratic forms
on functions $u\in C_c^\infty(B(x_0,r))$):
\begin{align*}
&H_{a+\ta,V}\le
(1+\eps)H_{a,V}+n\left(1+\frac{1}{\eps}\right)\|\ta\|^2_\infty\\
&=(1+\eps)\left(H_{a,V}+\frac{n}{\eps}\|\ta\|^2_\infty\right)
\le (1+\eps)\left(H_{a,V}+\frac{n^3r^2C^2}{\eps}\right).
\end{align*}
It follows that
$$
\la(B(x_0,r);H_{a+\ta,V})\le
(1+\eps)\left(\la(B(x_0,r);H_{a,V})+\frac{n^3r^2C^2}{\eps}\right).
$$
Using the IMS localization formula \eqref{E:IMS-localization}, 
we easily conclude that  $H_{a,V}$ is semi-bounded below on $C_c^\infty(\R^n)$, 
hence essentially self-adjoint due to Theorem \ref{T:self}.

Now the discreteness of the spectrum for $H_{a,V}$ immediately follows 
from the localization Theorem \ref{T:localization}.~$\ecarre$

\begin{corollary}\label{C:bound-magn}
Let $a$ be a magnetic potential such that $a\in\Liploc(\R^n)$
and $B=da$ is also in $\Liploc(\R^n)$ and bounded. Let
also $V\in L^\infty_{\textit{loc}}(\R^n)$ be
a bounded below electric potential. Then
$H_{a,V}$has a discrete spectrum
if and only if $V$ satisfies the Molchanov condition
$(M_c)$ with some $c>0$.
\end{corollary}

\textbf{Proof.} Theorem \ref{T:bound-magn} implies in
our case that $H_{a,V}$ has a discrete spectrum
if and only if this is true for $H_{0,V}$. By the Molchanov theorem 
the later is equivalent to the fulfillment of $(M_c)$ (for $V$) 
with some $c>0$. $\ecarre$

\ms
\begin{corollary}\label{C:const-magn}
Let $a$ be a magnetic potential which corresponds to
a constant magnetic field. Let
also $V\in L^\infty_{\textit{loc}}(\R^n)$ be
a bounded below electric potential. Then
$H_{a,V}$ has a discrete spectrum
if and only if $V$ satisfies the Molchanov condition
$(M_c)$ with some $c>0$.
\end{corollary}

\section{Sufficient conditions}\label{S:sufficient}

\subsection{Case $n=2$}\label{SS:sufficient2}
We will start by demonstrating the uncertainty principle argument
in the proof of the following well known Lemma 
(see e.g.\ \cite{Avron-Herbst-Simon}):

\begin{lemma}\label{L:uncertainty2}
Assume that $n=2$. Then
\begin{equation}\label{E:uncertainty2}
H_{a,0}\ge B(x)\quad \hbox{and} \quad H_{a,0}\ge -B(x),
\end{equation}
where $B$ as usual denotes the magnetic field produced by
the magnetic potential $a$.
The inequalities (\ref{E:uncertainty2}) hold
in the sense of operator inequalities (i.e.\ inequalities of
the quadratic forms) on $C_c^\infty(\R^2)$.
\end{lemma}

\textbf{Proof.} We have $H_{a,0}=P_1^2+P_2^2$ with
$[P_1,P_2]=-iB$. (See notations in Sect. \ref{SS:notation}.) Hence
for any $u\in C_c^\infty(\R^2)$
\begin{eqnarray*}
&(Bu,u)\le |(Bu,u)|&=|((P_1P_2-P_2P_1)u,u)|=|(P_2u,P_1u)-(P_1u,P_2u)|\\
&&= |2\hbox{Im} (P_1u,P_2u)|\le 2\|P_1u\|\|P_2u\|\le
\|P_1u\|^2+\|P_2u\|^2\\
&&=(H_{a,0}u,u).
\end{eqnarray*}
The second inequality in (\ref{E:uncertainty2}) follows
from the first one (by changing enumeration of coordinates).
$\ecarre$


\ms
\textbf{Proof of Theorem \ref{T:AHS-2dim}.} The condition $|B(x)|\to\infty$ means that either
$B\to +\infty$ or $-B\to +\infty$. In any of these cases the
condition $(d)$ of Theorem \ref{T:localization} is satisfied.
$\ecarre$

\begin{corollary}\label{C:BV-estimate2}
For $n=2$ and any $\de\in [-1,1]$ we have
\begin{equation}\label{E:BVde}
H_{a,V}\ge V+\de B
\end{equation}
\end{corollary}

\textbf{Proof.} Using the decomposition
$H_{a,V}=H_{a,0}+V$ and Lemma \ref{L:uncertainty2} we obtain
\begin{equation}\label{E:BV-uncertainty2}
H_{a,V}\ge B+V\quad \hbox{and} \quad H_{a,V}\ge -B+V.
\end{equation}
Multiplying the first inequality by $\ka\in [0,1]$, 
the second by $1-\ka$ and adding we obtain \eqref{E:BVde}
after denoting $1-2\ka$ by $\de$.
$\ecarre$

\begin{corollary}\label{C:abs-BV2}
Assume that $n=2$ and there exists $\de\in [-1,1]$ such that 
\begin{equation}\label{E:BV-to-infty}
V(x)+\de B(x)\to +\infty
\quad \textit{as} \quad x\to\infty.
\end{equation}
Then $H_{a,V}$ is essentially self-adjoint, semi-bounded below and
has a discrete spectrum.
\end{corollary}

\textbf{Proof.} Using Corollary \ref{C:BV-estimate2} we see that
the condition $(d)$ of Theorem \ref{T:localization} is fulfilled.
$\square$

\begin{corollary}\label{C:la-estimate2}
Assume that $n=2$ and $\Om\subset\R^2$ is an open set. Then
for any $\de\in [-1,1]$ 
\begin{equation}\label{E:la-estimate2}
\la(\Om;H_{a,V})\ge\underset{x\in\Om}\inf\{V(x)+\de B(x)\}.
\notag
\end{equation}
\end{corollary}

\textbf{Proof.} We can restrict the estimates given
in Lemma \ref{L:uncertainty2} and Corollary \ref{C:BV-estimate2}
to the functions from $C_c^\infty(\Om)$
for any $\Om\subset \R^2$. $\ecarre$ 

\medskip
\textbf{Proof of Theorem \ref{T:2dim}.} The result immediately
follows from Lemma \ref{L:uncertainty2} and
Theorem \ref{T:localization-cap} (with $\La(x)\equiv B(x)$).
$\ecarre$

\ms
{\textbf{Remark.}} V.~Ivrii noticed that if we take
any $\de\not\in [-1,1]$ then the result of 
the Corollary \ref{C:abs-BV2} does not hold any more.
Without loss of generality we can assume that 
$\de>1$. More precisely, we have

\begin{proposition}\label{P:Ivrii} {\rm(\rm V.~Ivrii)} There exists
a magnetic Schr\"odinger operator $H_{a,V}$ with $C^\infty$ potentials 
$a,V$ in $\R^2$ such that
$H_{a,V}\ge 0$, the spectrum of $H_{a,V}$ is not discrete,
but for any $\de>1$ we have $V(x)+\de B(x)\to +\infty$ as $x\to\infty$.
\end{proposition}

{\bf Proof.}
It is well known from the Landau calculation that
the bottom of the spectrum of the operator $H_{a,0}$ with constant 
magnetic field $B$ is precisely $|B|$ (see e.g.\ Sect.6.1.1. in \cite{Ivrii}
for a more general calculation). 

Let us take a sequence of constants $B_j$, $j=1,2,\dots$, $B_j>0$, $B_j\to +\infty$
as $j\to\infty$. Then define $V_j=-B_j$ and consider the Schr\"odinger operators
$H_j=H_{a_j,V_j}$, where $a_j$ is a magnetic potential corresponding to the constant
magnetic field $B_j$. Then the bottom of the spectrum of $H_j$ in $L^2(\R^2)$
is $0$ for all $j$. 

Note that the bottom of the spectrum of $H_j$ in $L^2(\R^2)$ can be defined as the bottom 
of the Dirichlet spectrum (see \eqref{E:lambda}) which involves only functions with 
a compact support. Therefore for any $\eps>0$ and any $x\in\R^2$ we can find $R_j>0$ 
such that $\la(B(x,R_j);H_j)<\eps$. In fact this does not depend either on the choice of 
the potential $a_j$, or on the center point $x$ due to the gauge invariance. 

Let us fix $\eps>0$ (e.g.\ take $\eps=1$).

Let us choose a sequence of points $x_j$, $j=1,2,\dots$, such that the balls $B(x_j,R_j+1)$
are  disjoint. Let us construct a function $B\in C^\infty(\R^2)$ such that $B(x)=B_j$ on
$B(x_j,R_j)$ and $B(x)\to +\infty$ as $x\to\infty$. By the Poincar\'e Lemma we can find 
a magnetic potential $a=a_kdx^k$ in $\R^2$ with $a_k\in C^\infty(\R^2)$, $k=1,2$, 
such that the corresponding magnetic field is $B(x)dx^1\wedge dx^2$ (or $B(x)$).
Define also $V(x)=-B(x)$ and consider the magnetic Schr\"odinger operator $H_{a,V}$
with $a$ and $V$ as constructed above.

It follows from Corollary \ref{C:BV-estimate2} (with $\de=1$) that $H_{a,V}\ge 0$.
On the other hand 
it is clear from the variational principle for the Dirichlet spectrum that 
there are infinitely many points of the spectrum of $H_{a,V}$ below $2\eps$.
Therefore the spectrum of $H_{a,V}$ is not discrete. At the same time
for any fixed $\de>1$ we have $V(x)+\de B(x)=(\de-1)B(x)\to +\infty$ as $x\to\infty$. 
$\square$

\subsection{Case $n\ge 3$} \label{SS:sufficient3}

The case $n\ge 3$ is substantially more complicated than the case
$n=2$, in particular because no growth condition on $|B|$ would
suffice for a reasonable estimate below for $H_{a,0}$ as was shown
by A.~Iwatsuka \cite{Iwatsuka1}. We will establish however that
appropriate regularity conditions on $B$, such as the ones imposed in
\cite{Avron-Herbst-Simon}, \cite{Dufresnoy}, \cite{Iwatsuka1},
can be incorporated in
the growth conditions for suitable effective potentials, so that
results similar to Theorem \ref{T:2dim} hold.

\subsubsection{The Iwatsuka identity}\label{SSS:Iwatsuka-id}
We will start with the following Lemma 
((6.2) on page 370 in \cite{Iwatsuka1}):

\begin{lemma}{\sc{(Iwatsuka identity)}}\label{L:Iwatsuka-id}
Assume that we are given an open set $\Om\subset\R^n$,
a  magnetic potential $a$
$($with components from $\Liploc(\Om)$$)$ 
and a set of real-valued functions
$A_{jk}\in \Liploc(\Om)$, $j,k=1,\dots,n$, $A_{kj}=-A_{jk}$. Then
\begin{equation}\label{E:Iwatsuka-id}
2\sum_{k<j}{\rm{Im}\,}(A_{kj}P_ku,P_ju)=
\left(\left[\sum_{k<j}A_{kj}B_{kj}\right]u,u\right)
+\sum_{k,j}\left(\frac{\pa A_{kj}}{\pa x^k}P_ju,u\right),
\end{equation}
for any $u\in W^{1,2}_{\textit{comp}}(\Om)$, i.e.\ $u\in L^2(\Om)$
such that $u$ has a compact support in $\Om$ and
$\nabla u\in (L^2(\Om))^n$ $($in the sense of
distributions$)$.
\end{lemma}

\textbf{Proof.} We reproduce the proof for the sake of
completeness and also because we have a different sign
convention in the definition of $a$ compared with
\cite{Iwatsuka1}.

An obvious approximation argument shows that it is sufficient
to consider $u\in C_c^\infty(\R^n)$. Then using integration by
parts and the commutation relations (\ref{E:commut}) we obtain
\begin{eqnarray*}
&&2\sum_{k<j}{\rm{Im}}\,(A_{kj}P_ku,P_ju)\\
&=&\frac{1}{i}\sum_{k<j}\left\{(A_{kj}P_ku,P_ju)-(A_{kj}P_ju,P_ku)\right\}\\
&=&\frac{1}{i}\sum_{k<j}\left((P_jA_{kj}P_k-P_kA_{kj}P_j)u,u)\right)\\
&=&\frac{1}{i}\sum_{k<j}\left\{(A_{kj}[P_j,P_k]u,u)+
\frac{1}{i}\left(\frac{\pa A_{kj}}{\pa x^j}P_ku,u\right)
-\frac{1}{i}\left(\frac{\pa A_{kj}}{\pa x^k}P_ju,u\right)\right\}\\
&=&\left(\left[\sum_{k<j}A_{kj}B_{kj}\right]u,u\right)
+\sum_{k,j}\left(\frac{\pa A_{kj}}{\pa x^k}P_ju,u\right). \quad
\ecarre
\end{eqnarray*}

\bs
It is natural to consider the functions $A_{jk}=-A_{kj}$ as
coefficients of a 2-form which we will call a \textit{dual field}
or a \textit{dual form}.

\ms
The Iwatsuka identity (\ref{E:Iwatsuka-id}) plays a very
important role in the arguments below. As A.~Iwatsuka noticed
in \cite{Iwatsuka1}, by choosing different
dual fields $A_{jk}$ we can obtain different estimates leading to
various sufficient conditions for the discreteness
of spectrum of the operator $H_{a,0}$. We will develop
this idea further by incorporating the electric potential
into the picture. In this way different effective potentials
$V_{\textit{eff}}$ emerge such that the Molchanov condition
for $V_{\textit{eff}}$ implies that $H_{a,V}$ is essentially
self-adjoint, semi-bounded below and has a discrete spectrum.
We will see that the  results of J.~Avron, I.~Herbst and B.~Simon
\cite{Avron-Herbst-Simon}, A.~Dufresnoy \cite{Dufresnoy}
and A.~Iwatsuka \cite{Iwatsuka1} about the discreteness
of spectrum can be extended to the case when both electric
and magnetic fields contribute to the localization
of the quantum particle.

\subsubsection{First choice}\label{SSS:Dufresnoy}

The first choice of the dual field which we will discuss, is
related to the Dufresnoy sufficiency result \cite{Dufresnoy}
and will lead to its generalization which takes into account
the electric field. In the situation of \cite{Dufresnoy}, considering
the operator $H_{a,0}$, we can take
\begin{equation}\label{E:Duf-choice}
A_{jk}(x)=\frac{B_{jk}(x)}{|B(x)|},
\end{equation}
where $|B|$ is defined as in (\ref{E:abs-B}).
(Note, however, that A.~Dufresnoy used different arguments in
\cite{Dufresnoy}.)
To be able to use
this choice we need to know that $B(x)\ne 0$ for large $|x|$ (i.e.\
if $|x|\ge R$ where $R>0$ is sufficiently large. This was
irrelevant in the situation when $|B(x)|\to\infty$ as $x\to\infty$
and $V=0$ which was considered in \cite{Dufresnoy}. In a more
general situation which will be considered here we do not want to
impose any a priori growth or non-vanishing requirement on $B$.
By this reason we will smooth down the field (\ref{E:Duf-choice}) at the 
places where $|B|$ is small and use the dual field given by 
\eqref{E:Duf-mod}.

\medskip
We will use the following notation:
$$
|\nabla B|=\left(\sum_{j<k}|\nabla B_{jk}|^2\right)^{1/2}.
$$
Similarly we define $|\nabla \beta|$ and $|\nabla A|$ below.

\medskip
Denote also $\beta_{jk}=B_{jk}/|B|$.
Differentiating $\beta_{jk}$ gives
\begin{equation}\label{E:nabla-beta-eq}
\nabla\beta_{jk}=|B|^{-1}\nabla B_{jk}-|B|^{-2}B_{jk}\nabla|B|.
\end{equation}
Using the inequality $|\nabla |B||\le |\nabla B|$, we obtain
\begin{equation}\label{E:nabla-beta-jk}
|\nabla\beta_{jk}|\le
|B|^{-1}|\nabla B_{jk}|+|B|^{-2}|B_{jk}||\nabla B|,
\end{equation}
hence by the triangle inequality
\begin{equation}\label{E:nabla-beta}
|\nabla\beta|\le 2|B|^{-1}|\nabla B|.
\end{equation}

Differentiating (\ref{E:Duf-mod}) we obtain
\begin{equation}\label{E:nabla-A-jk}
\nabla A_{jk}=
\chi(|B|)\nabla\beta_{jk}+\beta_{jk}\chi'(|B|)\nabla |B|,
\end{equation}
hence
$$
|\nabla A_{jk}|\le
\chi(|B|)|\nabla\beta_{jk}|+\chi'(|B|)|\beta_{jk}||\nabla B|.
$$
Therefore for $1/2\le |B|\le 1$ we have
\begin{equation}\label{E:nabla-A-jk-in}
|\nabla A_{jk}|\le |\nabla \beta_{jk}|+2|\beta_{jk}||\nabla B|
\end{equation}
and
\begin{equation}\label{E:nabla-A-in}
|\nabla A|\le |\nabla\beta|+2|\nabla B|
\le 2|B|^{-1}|\nabla B|+2|\nabla B|
\le 6|\nabla B|.
\end{equation}
Since $A=0$ for $|B|\le 1/2$, the estimate
(\ref{E:nabla-A-in}) holds if $|B|\le 1$.

Now introducing a majorant function
\begin{equation}\label{E:mixed-majorant}
M_B(x)=
\begin{cases}
|\nabla\beta_{jk}| &\text{if $|B|\ge 1$,}\\
6|\nabla B| &\text{if $|B|<1$},
\end{cases}
\end{equation}
 we see that
\begin{equation}\label{E:nabla-A-est-M-B}
|\nabla A(x)|\le M_B(x)\quad \textit{for all}\quad x\in\R^n.
\end{equation}
Note also that
$$
|B|^{-1}\le 2(1+|B|)^{-1} \quad\text{if $|B|\ge 1$,}
$$
and
$$
|\nabla B|\le 2(1+|B|)^{-1}|\nabla B| \quad\text{if $|B|\le 1$}.
$$
These estimates together imply that
$$
M_B(x)\le 12(1+|B|)^{-1}|\nabla B|
$$
for all $x$, hence (\ref{E:nabla-A-est-M-B}) gives
\begin{equation}\label{E:nabla-A-comb}
|\nabla A|\le 12(1+|B|)^{-1}|\nabla B|
\end{equation}
for all $x\in\R^n$.

\medskip
Now we are ready for 

\ms
\textbf{Proof of Theorem \ref{T:Duf-gen-maj}.} Let us use the Iwatsuka identity
(\ref{E:Iwatsuka-id}) with $A_{kj}$ given by (\ref{E:Duf-mod}).
Clearly
$$
\sum_{k<j}A_{kj}B_{kj}=\chi(|B|)|B|,
$$
and the first term in the right hand side of (\ref{E:Iwatsuka-id})
becomes $(\chi(|B|)|B|u,u)$. Since
$$
\chi(|B|)|B|\ge |B|-1, \quad \textit{and} \quad {|A_{jk}|\le 1},
$$
we obtain, using (\ref{E:A-majorant}), that
for any $u\in\Lip_c(\R)$:
\begin{eqnarray*}
&(|B|u,u)&\le 2\sum_{k<j}(|A_{kj}||P_ku|,|P_ju|)+
\sum_{k,j}\left(\left|\frac{\pa A_{kj}}{\pa x^k}\right||P_ju|,|u|\right)
+(u,u)\\
&&\le (n-1)\sum_{k=1}^n \|P_ku\|^2+
\sum_{j=1}^n (|P_ju|,X|u|)+(u,u)\\
&&=(n-1)h_{a,0}(u,u)+
\sum_{j=1}^n (|P_ju|,X|u|)+(u,u).
\end{eqnarray*}
Choosing an arbitrary $\eps>0$, we see that the middle term
in the right hand side is estimated by
$$
\eps\sum_{j=1}^n\|P_j\|^2+\frac{n}{4\eps}(X^2u,u)=
\eps h_{a,0}(u,u)+\frac{n}{4\eps}(X^2u,u).
$$
Therefore we obtain
$$
(|B|u,u)\le (n-1+\eps)h_{a,0}(u,u)+\frac{n}{4\eps}(X^2u,u)+(u,u),
$$
and
$$
h_{a,0}(u,u)\ge \frac{1}{n-1+\eps}(|B|u,u)-
\frac{n}{4\eps(n-1+\eps)}(X^2u,u)-\frac{1}{n-1+\eps}(u,u).
$$
Now we can apply Theorem \ref{T:localization-cap} with
$$
\La(x)=\frac{1}{n-1+\eps}|B|-
\frac{n}{4\eps(n-1+\eps)}X^2-\frac{1}{n-1+\eps},
$$
which immediately leads to the desired result. $\ecarre$

\medskip
Taking a specific majorant $X(x)$ in Theorem \ref{T:Duf-gen-maj}
we can make the result more specific. In this way we obtain
for example the following

\begin{corollary}\label{C:Duf-con-maj}
The result of Theorem \ref{T:Duf-gen-maj} holds if we replace
the majorant function $X$ in the definition of the effective
potential (\ref{E:Veff-de-gen}) either by $\sqrt{n-1}M_B$
or by $12\sqrt{n-1}(1+|B|^{-1})|\nabla B|$.
\end{corollary}

\textbf{Proof.} It is sufficient to notice that
$$
\sum_{k=1}^n
\left|\frac{\pa A_{kj}}{\pa x^k}\right|\le \sqrt{n-1}|\nabla A|,
$$
due to the Cauchy--Schwarz inequality. $\ecarre$

\medskip
We will need a notation for  domination between
two real-valued functions
$f,g$ on $S\subset\R^n$. Namely, we will write
\begin{equation}\label{E:dom-not}
f\prec g \qquad \textrm{or} \qquad g\succ f
\end{equation}
on $S$ if for any $\eps>0$ there exists $C(\eps)>0$ such that
\begin{equation}\label{E:dom}
f(x)\le\eps g(x)+C(\eps) \quad \textrm{for all}\quad x\in S.
\end{equation}
It is easy to see that $f\prec g$ and $g\prec h$ imply $f\prec h$.
If $f(x)\ge 0$, $g(x)>0$ for all $x\in S$,
$f,g$ are locally bounded and $g(x)\to +\infty$ as $x\to\infty$,
then the relation \eqref{E:dom-not} is equivalent to
$$
f(x)=o(g(x))\quad \textrm{as}\quad x\to\infty,
$$
i.e.\ $f(x)/g(x)\to 0$ as $x\to\infty$.

Note also that if $g$ is semi-bounded below on $S$, then
\begin{equation}\label{E:semi-dom}
f\prec g\quad \Longleftrightarrow\quad f\prec |g|
\quad\Longleftrightarrow\quad f\prec 1+|g|.
\end{equation}
The proof immediately follows from the implication
$$
g\ge -C\quad \Longrightarrow\quad |g|\le g+2C.
$$

\ms
Now we are ready for the formulation of an important corollary
of Theorem~\ref{T:Duf-gen-maj}. Note that though this theorem
does not explicitly
include any regularity requirements on $B$, they are implicit
in the requirements on the effective potential.
We will make the result more explicit (though weaker)
by invoking some explicit domination requirements
on the majorant $X$.

\begin{corollary}\label{C:X-cond}
Let us assume that there exists $\de\in [0,1)$ such that
the effective potential
\begin{equation}\label{E:V-eff-de}
\Veffde(x)=V(x)+\frac{\de}{n-1}|B(x)|
\end{equation}
is semi-bounded below and satisfies the Molchanov condition $(M_c)$
with some $c>0$ (in particular, this holds if $\Veffde(x)\to +\infty$
as $x\to\infty$).

In addition assume that the square of the majorant function
$X(x)$ from (\ref{E:A-majorant}) is dominated by $\Veffde$:
\begin{equation}\label{E:X-dom}
X^2\prec \Veffde \quad\textrm{on}\quad \R^n.
\end{equation}

\ms
Then the operator $H_{a,V}$ is essentially self-adjoint,
semi-bounded below and has a discrete spectrum.
\end{corollary}

\textbf{Remark.} Since $\Veffde$ is bounded below and locally bounded,
replacing $\Veffde$ by $1+|\Veffde|$ in \eqref{E:X-dom} leads to
an equivalent relation.

\ms
\textbf{Proof of Corollary \ref{C:X-cond}.} Let us choose
an arbitrary $\eps>0$. The condition (\ref{E:X-dom}) means that
for any $\kappa>0$ there exists $C(\kappa)>0$ such that
\begin{equation}\label{E:X2-Veff-est}
X^2(x)\le \kappa\Veffde(x)+C(\kappa), \quad x\in \R^n.
\end{equation}
Now for any $\eps>0$ and $\de\in [0,1)$ we can find $\eps>0$
and $\de'\in[0,1)$ such that
$$
\frac{\de'}{n-1+\eps}=\frac{\de}{n-1},
$$
so that
$$
\Veff^{(\de',\eps)}(x)=\Veffde(x)-
\frac{n\de'}{4\eps(n-1+\eps)}X^2(x)
$$
for $\Veffdeeps$ as in (\ref{E:Veff-de-gen}).
It follows from (\ref{E:X2-Veff-est}) that
$$
\Veff^{(\de',\eps)}(x)\ge (1-\kappa)\Veffde-C(\de,\eps,\kappa),
$$
hence $\Veff^{(\de',\eps)}$ is also semi-bounded below and  satisfies
$(M_c)$ with the same $c>0$ as for $\Veffde$. Therefore the desired
result follows from Theorem \ref{T:Duf-gen-maj}. $\ecarre$

\medskip
Now we will replace the conditions
on the majorant $X$ in Corollary \ref{C:X-cond} by an
explicit choice of $X$ which is
of A.~Dufresnoy type \cite{Dufresnoy}.
In this way we get a theorem which improves
the sufficiency result of \cite{Dufresnoy} adding electric field
and capacity to the picture.

\begin{theorem}\label{T:o-Veff}
Let us assume that there exists $\de\in [0,1)$ such that
the effective potential (\ref{E:V-eff-de})
is semi-bounded below and satisfies the Molchanov condition $(M_c)$
with some $c>0$ (in particular, this holds if $\Veffde(x)\to +\infty$
as $x\to\infty$).

In addition assume that
one of the following conditions  $(a)$, $(b)$ is satisfied:

\ms
$($a$)$ $M_B^2\prec\Veffde$,

\ms
or, in other words,

$
|\nabla\beta|^2\prec \Veffde\quad
\textrm{on}\quad \{x|\,|B(x)|\ge 1\},
$

and

$
|\nabla B|^2\prec \Veffde(x)\quad
\textrm{on}\quad \{x|\,|B(x)|<1\};
$

\ms
$($b$)$ $(1+|B|)^{-2}|\nabla B|^2\prec \Veffde$ on $\R^n$.

\ms
Then the operator $H_{a,V}$ is essentially self-adjoint,
semi-bounded below and has a discrete spectrum.
\end{theorem}

\ms
\textbf{Proof.} According to
(\ref{E:nabla-A-est-M-B}) and (\ref{E:nabla-A-comb})
the conditions $(a)$ or $(b)$ imply that
we can take the majorant
$$
X(x)=M_B(x)\quad \textit{or}\quad 12(1+|B(x)|)^{-1}|\nabla B(x)|
$$
respectively, and apply Corollary \ref{C:X-cond}. $\ecarre$

The following even more explicit result also
improves the sufficiency result by A.~Dufresnoy \cite{Dufresnoy}.

\begin{theorem}\label{T:mod-Duf}
Let us assume that $B_{jk}\in \Lip(\R^n)$ for all $j,k$,
and the following conditions are satisfied:
\begin{equation}\label{E:nabla-B-bound}
|\nabla B|\le C \quad\textit{if}\quad |B|\le 1;
\end{equation}
\begin{equation}\label{E:Duf-beta-est}
|\nabla\beta|=o(|B|^{1/2}) \quad\textit{as}\quad |B|\to\infty.
\end{equation}
Assume also that there exists $\de\in[0,1)$
such that the effective potential (\ref{E:V-eff-de}) is
semi-bounded below and satisfies the Molchanov condition $(M_c)$
with some $c>0$. Then the operator $H_{a,V}$ is essentially
self-adjoint, semi-bounded below and has a discrete spectrum.
\end{theorem}

\textbf{Proof.} The result easily follows from
Theorem \ref{T:o-Veff} if $V\ge 0$ (or if $V$ is semi-bounded
below). In the general case it can be deduced from
Theorem \ref{T:Duf-gen-maj} or Corollary \ref{C:Duf-con-maj}
by the same argument which was used in the proof of
Corollary \ref{C:X-cond}, except that in this case
instead of (\ref{E:X2-Veff-est}) we should estimate
$X^2(x)$ by $\kappa |B(x)|+C(\kappa)$.  $\ecarre$

\bigskip
\textbf{Proof of Theorem \ref{T:B-and-V}.}
To deduce this theorem from Theorem \ref{T:mod-Duf}
it is sufficient to show that the condition
(\ref{E:Duf-beta-est}) imposed on $\beta$ in
Theorem \ref{T:mod-Duf}, follows from
the condition $(B_{3/2}^0)$ which is imposed on $B$
and means that
$$
|\nabla B|=o\left((1+|B|)^{3/2}\right)\quad
\textit{as}\quad |B|\to\infty.
$$
This implication immediately follows from
the estimate (\ref{E:nabla-beta}). $\ecarre$

\bigskip
\textbf{Remark 1.} It is also possible to prove Theorem \ref{T:B-and-V}
and all previous theorems of this subsection
(except the ones with conditions which explicitly include
$\beta_{jk}$) by the following choice
of the dual field:
$$
A_{jk}(x)=\frac{B_{jk}(x)}{\lBxr},
$$
where $\lBr=(1+|B|^2)^{1/2}$. This choice leads to
the arguments and estimates which are similar to the ones used
above in the proof of Theorem \ref{T:Duf-gen-maj}.

\medskip
\textbf{Remark 2.} Sufficient measure conditions similar to
Corollaries \ref{C:3dim-mes} and \ref{C:3dim-elem} hold
for all types of effective potentials discussed above.

\subsubsection{Second choice}\label{SSS:Iwatsuka}
Now we will discuss the choice of the dual field $A_{jk}$
which is associated with the field
suggested by A.~Iwatsuka \cite{Iwatsuka1}:
$$
A_{jk}(x)=\frac{B_{jk}(x)}{|B(x)|^2}.
$$
It is proved in \cite{Iwatsuka1} that this choice leads to
a weakest regularity condition on $B$ which guarantees the
discreteness of spectrum for $H_{a,0}$ provided
$|B(x)|\to\infty$ as $x\to\infty$.

We will improve the result of \cite{Iwatsuka1} by adding a term
taking into account the electric field.
Hence by the same reason as above we will slightly modify
this field as follows:
\begin{equation}\label{E:Ajk-Iwatsuka}
A_{jk}(x)=\frac{B_{jk}(x)}{\lBxr^2}.
\end{equation}
We will assume that $B_{jk}\in \Liploc(\R^n)$ for all $j,k$.
A.~Iwatsuka \cite{Iwatsuka1} assumes also that
the condition $(B_2^0)$ from Sect.\ref{SS:suff3}
holds.
Choosing an arbitrarily small $r>0$ denote
\begin{equation}\label{E:Iw-eps}
\eps_x=\sup_{y\in B(x,r)}\frac{1+|\nabla B(y)|}
{\langle B(y)\rangle^2}.
\end{equation}
The conditions \eqref{E:B-to-infty} and  $(B_2^0)$ together are
equivalent to the relation
\begin{equation}\label{E:eps-to-0}
\eps_x\to 0\quad \textit{as}\quad x\to\infty.
\end{equation}
It is proved in \cite{Iwatsuka1} that this implies that $H_{a,0}$
has a discrete spectrum.

We will not a priori require either (\ref{E:B-to-infty}) or $(B_2^0)$
to be satisfied, but we will include $\eps_x$ above into an
effective potential so that possible violation of
(\ref{E:eps-to-0}) is compensated by the electric field.
More precisely, we will prove

\begin{theorem}\label{T:Iwa-mod}
Let us assume that $B_{jk}\in \Lip(\R^n)$ for all $j,k$ and
define an effective potential
\begin{equation}\label{E:Veff-Iw}
\Veffde(x)=V(x)+\frac{\de}{n-1}\, \eps_x^{-1/2}(1-\eps_x).
\end{equation}
If there exists $\de\in [0,1)$ such that
$\Veffde$ is semi-bounded below
and satisfies the Molchanov condition $(M_c)$ with a positive
$c>0$, then $H_{a,V}$ is essentially self-adjoint, semi-bounded
below and has a discrete spectrum.
\end{theorem}

\textbf{Proof.} Note first that due to the triangle inequality
$$
|B|\le\lBr \le 1+|B|,
$$
so $\lBr$ can be replaced by $|B|$ (and vice versa)
in the domination relations.

We will again use the Iwatsuka identity \eqref{E:Iwatsuka-id}.
With the choice \eqref{E:Ajk-Iwatsuka} we have
\begin{equation}\label{E:AkjBkj}
\sum_{k<j}A_{kj}B_{kj}=\frac{|B|^2}{\lBr^2}=\frac{\lBr^2-1}{\lBr^2}
=1-\frac{1}{\lBr^2}.
\end{equation}

Let us fix $r>0$ and choose $x\in\R^n$. Then using the inequality
$$
|A_{kj}(y)|\le \lByr^{-1}\le \eps_x^{1/2},\quad y\in B(x,r),
$$
we obtain for any $u\in C_c^\infty(B(x,r))$
for the left hand side of \eqref{E:Iwatsuka-id}:
\begin{align}\label{E:Iw-left-est}
&2\underset{k<j}{\sum}{\textrm{Im}\,}(A_{kj}P_ku,P_ju)
\le 2\underset{k<j}{\sum}(|A_{kj}||P_ku|,|P_ju|)\\
&\le 2\eps_x^{1/2} \underset{k<j}\sum(|P_ku|,|P_ju|)
\le (n-1)\eps_x^{1/2} h_{a,0}(u,u).\notag
\end{align}

Now let us estimate the last term in \eqref{E:Iwatsuka-id}.
Using the inequality $|\nabla |B||\le |\nabla B|$, we obtain
\begin{eqnarray*}
&|\nabla A_{kj}|&=\left|\lBr^{-2}\nabla B_{kj}-
2B_{kj}\lBr^{-4}|B|\nabla |B|\right|
\notag\\
&&\le \lBr^{-2}|\nabla B_{kj}|+
2|B_{kj}|\lBr^{-4}|B||\nabla B|
\notag\\
&&\le \lBr^{-2}|\nabla B_{kj}|+
2|B_{kj}|\lBr^{-3}|\nabla B|.
\notag
\end{eqnarray*}
Therefore by the triangle inequality
$$
|\nabla A|\le 3\lBr^{-2}|\nabla B|.
$$
Now using the Cauchy-Schwarz inequality we find:
$$
\sum_k\left|\frac{\pa A_{kj}}{\pa x^k}\right|\le
\sqrt{n-1} |\nabla A|\le 3\sqrt{n-1}\lBr^{-2}|\nabla B|.
$$
Hence we can estimate the last term in \eqref{E:Iwatsuka-id}
as follows:
\begin{align}
&\left|\sum_{k,j}\left(\frac{\pa A_{kj}}{\pa x^k}P_ju,u\right)\right|
\le 3\eps_x\sqrt{n-1} \sum_j(|P_ju|,|u|)
\notag\\
&\le 3\eps_x\sqrt{n-1}
\left(\kappa_x\sum_j
\|P_ju\|^2+\frac{n}{4\kappa_x}\|u\|^2\right)
\notag\\
&=3\eps_x\sqrt{n-1}\left(\kappa_x h_{a,0}(u,u)+
\frac{n}{4\kappa_x}\|u\|^2\right),
\notag
\end{align}
where $u\in C_c^\infty(B(x,r))$ and $\kappa_x>0$ is arbitrary.

Let us fix $\eps>0$ and then choose $\kappa_x$ so that
$$
3\eps_x\sqrt{n-1}\kappa_x=\eps \eps_x^{1/2}.
$$
This leads to the estimate
\begin{equation}\label{E:Iw-right2-est}
\left|\sum_{k,j}\left(\frac{\pa A_{kj}}{\pa x^k}P_ju,u\right)\right|
\le
\eps\eps_x^{1/2}h_{a,0}(u,u)+\frac{9n(n-1)}{4\eps}\eps_x^{3/2}\|u\|^2.
\end{equation}
Now from the identities \eqref{E:Iwatsuka-id}, \eqref{E:AkjBkj} and from
the estimates \eqref{E:Iw-left-est}, \eqref{E:Iw-right2-est} we obtain
for any $u\in C_c^\infty(B(x,r))$
\begin{equation*}
(n-1+\eps)h_{a,0}(u,u)\ge
\left[(1-\eps_x)\eps_x^{-1/2}-\frac{9n(n-1)\eps_x}{4\eps}\right](u,u).
\end{equation*}
(Here we used the obvious estimate $\lByr^{-2}\le \eps_x$ if
$y\in B(x,r)$.)
The last estimate gives a lower bound for the bottom
of the Dirichlet spectrum of $H_{a,0}$ on any ball $B(x,r')$
with $r'<r$:
\begin{equation*}
\la(B(x,r');H_{a,0})\ge
(n-1+\eps)^{-1}\left[\eps_x^{-1/2}(1-\eps_x)-\frac{9n(n-1)\eps_x}{4\eps}\right].
\end{equation*}
If $\eps_x<1$, we can replace $\eps_x$ by 1 in the last term to get
\begin{equation}\label{E:Dir-above}
\la(B(x,r');H_{a,0})\ge
(n-1+\eps)^{-1}\left[\eps_x^{-1/2}(1-\eps_x)-\frac{9n(n-1)}{4\eps}\right].
\end{equation}
This also obviously holds if $\eps_x\ge 1$.

Now we can apply the IMS localization formula argument (see proof of
Theorem \ref{T:localization}) by considering a finite multiplicity
covering of $\R^n$ by balls of radius $r/2$, to conclude that
$H_{a,0}\ge \La(x)$ with
$$
\La(x)=(n-1+\eps)^{-1}\eps_x^{-1/2}(1-\eps_x)-C,
$$
where $C=C(r,n,\eps)$. Then the statement of the Theorem follows from
Theorem \ref{T:localization-cap}. $\ecarre$

\ms
\begin{corollary}\label{C:Iwa-mod}
Let us assume that $B$ satisfies the Iwatsuka conditions
\eqref{E:B-to-infty} and $(B_2^0)$
(or, equivalently, $\eps_x\to 0$ as $x\to\infty$).
Then the statement of the Theorem \ref{T:Iwa-mod} holds with
the effective potential
$$
\Veffde(x)=V(x)+\frac{\de}{n-1}\,\eps_x^{-1/2}.
$$
\end{corollary}

\textbf{Remark 1.}
If $V=0$ then this Corollary immediately implies the Iwatsuka
result \cite{Iwatsuka1} about the sufficient condition of
the discreteness of the spectrum because then the Molchanov
condition $(M_c)$ for the function $x\mapsto\eps_x^{-1/2}$
is fulfilled automatically.

\ms
\textbf{Remark 2.} If $|B(x)|\to\infty$ as $x\to\infty$ and $|B(x)|$
varies sufficiently slowly  (e.g.\ if $|\nabla |B||\prec |B|$,
which in turn holds e.g.\ if $B$ satisfies the condition $(B_1^0)$), then
$\eps_x^{-1/2}$ becomes equivalent to $|B(x)|$, so the results
of Sect.\ref{SSS:Iwatsuka} agree with the results of
Sect.\ref{SSS:Dufresnoy}.

\subsubsection{Third choice}

Here we will discuss the choice of the dual field $A_{jk}$
which leads to a generalization of the result by
J.~Avron, I.~Herbst and B.~Simon \cite{Avron-Herbst-Simon}.
In this choice $A_{jk}$ are constants, though we argue on balls
of a fixed radius $r>0$ and these constants may
also depend on the choice of the ball.

The advantage of the result obtained in this way is that
no local smoothness of $B$ is required. Even the continuity of $B$
is not needed, though we will still maintain the requirement
$a_j\in\Liploc(\R^n)$, hence
$B_{jk}\in L^\infty_{\textit{loc}}(\R^n)$.

Let us choose a finite multiplicity covering of $\R^n$ by balls
$B(\ga,r)$, $\ga\in\Ga$. (For example $\Ga$ may be an appropriate
lattice in $\R^n$.) Then for any $\ga\in\Ga$ chose
$A_{kj}(\ga)=-A_{jk}(\ga)$ such that $|A(\ga)|=1$ and denote
\begin{equation}\label{E:r-ga}
r_\ga=\inf_{y\in B(\ga,r)}\sum_{k<j}A_{kj}(\ga)B_{kj}(y).
\end{equation}
We are interested to make the  numbers $r_\ga$ as big as possible.
Clearly
$$
r_\ga\le \sup_{y\in B(\ga,r)}|B(y)|,
$$
and this inequality is close to equality if $B$ is almost
constant in $B(\ga,r)$ and we have chosen
$A_{kj}(\ga)=B_{kj}(\ga)/|B_{kj}(\ga)|$.

It is proved in \cite{Avron-Herbst-Simon} that if we can choose
$A_{kj}(\ga)$ so that
\begin{equation}\label{E:r-to-infty}
r_\ga\to\infty \quad\textit{as}\quad \ga\to\infty,
\end{equation}
then the operator $H_{a,0}$ has a discrete spectrum.
We will improve this result by adding an electric field into the
picture. Denote
\begin{equation}\label{E:rx}
r(x)=\min \{r_\ga|\, |x-\ga|\le r\}.
\end{equation}

\begin{theorem}\label{T:AHS}
Denote
\begin{equation}\label{E:Veff-AHS}
\Veffde(x)=V(x)+\frac{\de}{n-1}\; r(x),
\end{equation}
and assume that there exists $\de\in [0,1)$ such that the
effective potential \eqref{E:Veff-AHS} is semi-bounded below and
satisfies the Molchanov condition $(M_c)$ with some $c>0$. Then
the operator $H_{a,V}$ is essentially self-adjoint, semi-bounded
below and has a discrete spectrum.
\end{theorem}

\textbf{Proof.} Let us use the Iwatsuka identity \eqref{E:Iwatsuka-id}
for $u\in C_c^\infty(B(\ga,r))$ and with $A_{kj}(x)=A_{kj}(\ga)$.
The last term in \eqref{E:Iwatsuka-id} vanishes. The first term in
the right hand side of \eqref{E:Iwatsuka-id} is estimated from
below by $r_\ga(u,u)$ and the left hand side is estimated from
above by $(n-1)h_{a,0}(u,u)$ (see arguments given in the
proof of Theorem \ref{T:Duf-gen-maj}). Therefore we obtain:
\begin{equation}\label{E:Iwa1-AHS}
h_{a,0}(u,u)\ge \frac{r_\ga}{n-1}(u,u),\quad u\in
C_c^\infty(B(\ga,r)).
\end{equation}
Now using the IMS localization formula as in the proof of
Theorem \ref{T:Iwa-mod} we conclude that
$$
H_{a,0}\ge \frac{r(x)}{n-1}-C,
$$
and it remains to apply Theorem \ref{T:localization-cap}.
$\ecarre$

\section{Operators on manifolds}\label{S:manifolds}

Let $(M,g)$ be a Riemannian manifold (i.e.\ $M$ is a $C^\infty$-manifold,
$(g_{jk})$ is a Riemannian metric on $M$), $\dim M=n$.
We will always assume that $M$ is connected. We will also assume
that we are given a {\it positive smooth measure} $d\mu$, i.e.\
a measure which has a $C^\infty$ positive density $\rho(x)$
with respect to the Lebesgue measure $dx=dx^1\dots dx^n$
in any local coordinates $x^1,\dots,x^n$, so we will write
$d\mu=\rho(x)dx$.
This measure may be completely independent of the Riemannian
metric, but may of course coincide with the canonical measure $d\mu_g$
induced by the metric (in this case $\rho=\sqrt{g}$ where
$g=\det(g_{jk})$, so locally $d\mu_g=\sqrt{g}dx$).

Denote $\La^p_{(k)}(M)$ the set of all $k$-smooth (i.e.\ of the class $C^k$)
complex-valued $p$-forms on $M$. We will write $\La^p(M)$ instead
of $\La^p_{(\infty)}(M)$.
A {\it magnetic potential} is a
real-valued 1-form $a\in \La^1_{(1)}(M)$. So in local coordinates
$x^1,\dots, x^n$ it can be written as in (\ref{E:a-form})
where $a_j=a_j(x)$ are real-valued $C^1$-functions of the
local coordinates.

The usual differential can be considered as a first order
differential operator
\begin{equation}\notag
d:C^\infty(M)\longrightarrow \La^1(M).
\end{equation}
The deformed differential (compare (\ref{E:mag-diff}))
\begin{equation}\notag
d_a:C^\infty(M)\longrightarrow \La^1_{(1)}(M),\qquad u\mapsto du+iua,
\end{equation}
is also well defined.

The Riemannian metric $(g_{jk})$ and the measure $d\mu$ induce an inner
product in the spaces of smooth forms with compact support
in a standard way.
The corresponding completions 
are Hilbert spaces which we will denote
$L^2(M)$  for functions and $L^2\La^1(M)$ for 1-forms.
These spaces depend on the choice of the metric $(g_{jk})$ and
the measure $d\mu$. However we will skip this dependence in the
notations of the spaces for simplicity of notations. 

The corresponding local spaces will be denoted $L^2_{loc}(M)$
and $L^2_{loc}\La^1(M)$ respectively. These spaces do not depend
on the metric or measure.

Formally adjoint operators to the differential operators
with sufficiently smooth coefficients are well defined
through the inner products above. In particular,
we have an operator
\begin{equation}\notag
d_a^*:\La^1_{(1)}(M)\longrightarrow C(M),
\end{equation}
defined by the identity
\begin{equation}
(d_au,\om)=(u,d_a^*\om),\ u\in C^\infty_c(M),\ \om\in\La^1_{(1)}(M).
\notag
\end{equation}
(Here $C_c^\infty(M)$ is the set of all
$C^\infty$ functions with compact support on $M$.)

Therefore we can define the magnetic Laplacian  $\De_a$
(with the potential $a$) by the formula
\begin{equation}\notag
-\De_a=d_a^* d_a: C^\infty(M)\longrightarrow C(M).
\end{equation}
Now the
\textit{magnetic Schr\"odinger operator} on $M$ is defined as
\begin{equation}\label{E:Hmag}
H=H_{a,V}=-\De_a+V,
\end{equation}
where $V\in L^\infty_{loc}(M)$, i.e.\ $V$ is a locally bounded measurable
function which is called {\it the electric potential}. We will always
assume $V$ to be real-valued. Then $H$ becomes a symmetric
operator in $L^2(M)$ if we consider it on the domain
$C_c^\infty(M)$.

Note that for $a=0$  the operator $\De_a$ becomes
a generalized Laplace-Beltrami
operator $\De$ on scalar functions on $M$ and it can be locally
written in the form
\begin{equation}\label{E:Laplace-Beltrami}
\De u=\frac{1}{\rho}\frac{\pa}{\pa x^j}(\rho g^{jk}
\frac{\pa u}{\pa x^k}).
\end{equation}

The following  expressions for $H_{a,V}$ are also useful
(\cite{Shubin3}):
\begin{equation}\label{E:HaV-explicit}
H_{a,V}u=-\De u -2i\langle a, du\rangle+(id^*a+|a|^2)u+Vu,
\end{equation}
and in local coordinates

\begin{equation}\label{E:HaV-local}
H_{a,V}u=-\frac{1}{\rho}\left(\frac{\partial}{\partial x^j}+ia_j\right)
\left[\rho g^{jk}
\left(\frac{\partial}{\partial x^k}+ia_k\right)u\right]+Vu,
\end{equation}
or in slightly different notations
\begin{equation}\label{E:HaV-local-2}
H_{a,V}u=
\frac{1}{\rho}P_j
[\rho g^{jk}P_ku]+Vu=
\frac{1}{\rho}(D_j+a_j)
[\rho g^{jk}(D_k+a_k)u]+Vu,
\end{equation}
where $D_j=-i\partial/\partial x_j$.

Now we need a condition on the Riemannian manifold $(M,g)$
which would allow extending the results above to a more general
context. This is the condition of \textit{bounded geometry}
which means that the injectivity radius $r_{inj}$ is positive
and all covariant derivatives of the curvature tensor $R$
are bounded:
\begin{eqnarray}\notag
&(a) &r_{inj}>0;\\
\notag
&(b) &|\nabla^m R|\le C_m; m=0,1,\dots.
\end{eqnarray}
Here the norm $|\cdot|$ of tensors $\nabla^m R$ is measured with
respect to the given Riemannian metric $g$. (See more details about
these conditions and their use in \cite{Kondrat'ev-Shubin}.)

We will also impose a bounded geometry condition on the measure
$d\mu$. We will say that \textit{the triple} $(M,g,d\mu)$ 
\textit{has bounded geometry} if $(M,g)$ is a manifold 
of bounded geometry, and for a small $r>0$ in local geodesic
coordinates in any ball $B(x,r)$ we have $d\mu=\rho(x)dx$
where $\rho\ge\eps>0$ and for any multiindex $\al$
we have $|\pa^\al\rho|\le C_\al$ with the constants 
$\eps$, $C_\al$ independent of $x$.
In particular this automatically holds for the Riemannian
measure $d\mu=\sqrt{g}dx$ if $(M,g)$ has bounded geometry.

The requirement on the measure was not needed for $\R^n$
which reflects the fact that our methods do not work
as well for manifolds as they do for $\R^n$.

Some conditions
at infinity are needed to guarantee that the operator $H_{a,V}$
is essentially self-adjoint in $L^2(M)$
- see e.g.\ M.~Shubin \cite{Shubin3, Shubin00} and references there for such conditions.
The essential self-adjointness result by A.~Iwatsuka
\cite{Iwatsuka2} can also be extended to manifolds
of bounded geometry. The result which is of particular importance
for us is essential self-adjointness of any semi-bounded below
magnetic Schr\"odinger operator on any complete
Riemannian manifold (see \cite{Shubin00}), in particular on any manifold
of bounded geometry.  

For a triple $(M,g,d\mu)$ of bounded geometry we can define a capacity
of a compact set $F\subset B(x,r)$ for a small $r>0$ by use of
geodesic coordinates or by use of norms induced by the metric $g$ and
the measure $d\mu$.
Using geodesic coordinates for different balls $B(x,r)$ and $B(x',r)$
to measure $\capa(F)$ for $F\subset B(x,r)\cap B(x',r)$ leads
to equivalent results, so it does not affect our formulations.
Using Riemannian norms of tensors in the definition of capacity
also leads to an equivalent result.

After these explanations all results formulated above make sense and can be
extended to the triples of bounded geometry with minor changes
(some constants depending on geometry appear in the formulations). In case 
when $a=0$ and with the Riemannian measure $d\mu$ this was done
by the authors in \cite{Kondrat'ev-Shubin}. We will give
examples of such extensions in Sect.\ref{S:manifolds}.

Let $(M,g, d\mu)$ be a triple
with bounded geometry. Let us chose a ball $B(x_0,r)$ of a fixed 
sufficiently small radius $r>0$, and let $x^1,\dots,x^n$ be local geodesic 
coordinates in this ball.

\begin{lemma}\label{L:ball-est} 
Under the conditions above there exists $C=C(M,g,d\mu)$ such that
\begin{equation}\label{E:ball-est}
C^{-1}\sum_{j=0}^n\|P_ju\|^2_0\le h_{a,0}(u,u)\le C\sum_{j=0}^n\|P_ju\|^2_0, 
\quad u\in C_c^\infty(B(x_0,r)),
\end{equation}
where $P_j$ are defined by \eqref{E:Pj} in the geodesic coordinates,
and  $\|\cdot\|_0$ means the norm in $L^2(B(x_0,r); dx^1\dots dx^n)$.
\end{lemma}

\textbf{Proof.} The result immediately follows from the local presentation
\eqref{E:HaV-local-2}  for $H_{a,V}$ and from the bounded geometry requirements. 
$\square$

This Lemma allows an easy extension of all local estimates in the balls 
of a fixed small radius $r>0$ to the case of the triples $(M,g,d\mu)$ of
bounded geometry. Then we need to use extensions of the localization results
from Sect.\ref{S:localization}. They can indeed be extended due to Gromov's
observation on the existence of good coverings \cite{Gromov}
and the manifolds extension of
the IMS localization formula \eqref{E:IMS-localization} -- see the details in
\cite{Shubin1, Shubin96, Kondrat'ev-Shubin}.

Now we can formulate the simplest result for 2-dimensional manifolds.

\begin{theorem}\label{T:2dim-M} Assume that we have 
a bounded geometry triple $(M,g,d\mu)$ with $n=\dim M=2$.
Then there exists $\de_0>0$ depending on the triple $(M,g,d\mu)$,
such that the following holds. 

Let $H_{a,V}$ be a magnetic
Schr\"odinger operator on $M$ 
and there exists $\de\in (-\de_0,\de_0)$ such that
the effective potential 
given by (\ref{E:Veff}) with $B$ identified
with the ratio $B/d\mu$,
is semi-bounded below and  satisfies
the Molchanov condition $(M_c)$
with some $c>0$. Then $H_{a,V}$ is essentially self-adjoint,
semi-bounded below and has a discrete spectrum.
\end{theorem}

\textbf{Proof.} The estimate $H_{a,0}\ge \de B$ follows from 
Lemmas \ref{L:ball-est} and \ref{L:uncertainty2}
on balls $B(x,r)$. Then the IMS-localization formula 
\eqref{E:IMS-localization} used for a 
finite multiplicity covering of $M$ by such balls (of the same radius)
leads to a global estimate of the form
\begin{equation*}
H_{a,0}\ge \de B-C,
\end{equation*}
where $C=C(M,g,d\mu)$. Now the desired result immediately follows
from the  manifold version of Theorem \ref{T:localization-cap}.
$\square$

Let us formulate an extension of Theorem \ref{T:Duf-gen-maj} to manifolds.

\ms
Let us define a smoothed direction of the magnetic field as a 
2-form (or a skew-symmetric (0,2)-tensor):
\begin{equation}\label{E:Duf-mod-M}
A(x)=\chi(|B(x)|)\frac{B(x)}{|B(x)|},
\end{equation}
where $\chi\in \Lip([0,\infty))$, $\chi(r)=0$ if $r\le 1/2$,
$\chi(r)=1$ if $r\ge 1$, $\chi(r)=2r-1$ if $1/2\le r\le 1$, so
$0\le \chi(r)\le 1$ and $|\chi'(r)|\le 2$
for all $r$. The norm $|B|$  is measured by 
the use of the given Riemannian metric $g$.

\begin{theorem}\label{T:Duf-gen-maj-M}
Let us assume that $B\in \Lip_{loc}(\R^n)$,
$A$ is defined by (\ref{E:Duf-mod-M}),
and a positive measurable function
$X(x)$ in $\R^n$ satisfies
\begin{equation}\label{E:A-majorant-M}
|\nabla A(x)|\le X(x), \quad x\in M,
\end{equation}
where $\nabla$ means the covariant derivative. 
Then there exist constants $\de_0>0$ and $C_0>0$ depending
only on the triple $(M,g,d\mu)$, such that the following is true.
If there exists $\de\in(-\de_0,\de_0)$ such that the effective potential  
\begin{equation}\label{E:Veff-de-gen-M}
\Veff(x)=V(x)+\de |B(x)|-C_0\de X^2(x)
\end{equation}
is bounded below and satisfies the Molchanov condition $(M_c)$
with some $c>0$ (in particular, this holds if
$\Veff(x)\to +\infty$ as $x\to\infty$),
then $H_{a,V}$ is essentially self-adjoint,
semi-bounded below and has a discrete spectrum.
\end{theorem}

In case $B=0$ and $V$ semi-bounded below the condition on 
$V$ in this theorem becomes necessary and sufficient due 
to the extension of Molchanov theorem given in \cite{Kondrat'ev-Shubin}.

Proof of Theorem \ref{T:Duf-gen-maj-M} is similar to the proof of
Theorem \ref{T:2dim-M}. Other results from Sect.\ref{S:sufficient}
have similar extensions as well.

\section{Other results}

A review of other results on the discreteness of spectrum of the Schr\"odinger operators
and related topics can be found in \cite{Kondrat'ev-Shubin}. Here we will restrict
ourselves to a few specific remarks concerning magnetic Schr\"odinger operators.

Cwickel-Lieb-Rozenblum (CLR) type estimates for the number of negative eigenvalues 
for magnetic Schr\"odinger operators have been proved by J.~Avron, I.~Herbst and B.~Simon 
\cite{Avron-Herbst-Simon} (see also \cite{Simon79b}), though with the right hand side
independent of the magnetic field. The proof used heat kernel estimates 
based on the Feynman--Kac formula as in the paper by E.~Lieb \cite{Lieb76}, and also
the diamagnetic inequality (see \cite{Simon79}, \cite{Hess77} or 
\cite{Cycon-Froese-Kirsch-Simon}, Sect.1.3). 
An analytic proof was provided by M.~Melgaard and G.~Rozenblum \cite{Melgaard-Rozenblum}. 

The CLR estimates used for $H_{a,V}-M$ for arbitrary $M>0$ obviously imply sufficient conditions
for the discreteness of spectrum (namely, the finiteness 
of the right hand sides of these estimates
for all $M$). However, the  above mentioned results still 
do not allow to take into account
any interaction between the electric and magnetic fields. 

Under some stronger conditions on the fields it is possible to obtain even asymptotics for
the counting function $N(\lambda;H_{a,V})$ for the  eigenvalues of $H_{a,V}$. One of the first results
of this kind is due to Y.~Colin de Verdi\`ere \cite{Colin-de-Verdiere}. Numerous further results
on such asymptotics can be found in (or extracted from) the book of V.~Ivrii \cite{Ivrii}
(see also references there).

The Lieb--Thirring inequalities give explicit estimates for sums of powers of the 
negative eigenvalues, or, in other words, for the $l^p$ norms of the sequence 
of these eigenvalues. If $p=\infty$ this means an estimate for the number
of negative eigenvalues as in the case of the CRL estimate. Under some conditions 
Lieb--Thirring type inequalities were obtained for $H_{a,V}$ by L.~Erd\"os \cite{Erdos}
and for a similar Pauli operator by A.V.~Sobolev \cite{Sobolev}.

A Feynman type estimate for $\Tr(\exp (-tH_{a,V}))$ in explicit purely classical
terms was obtained by J.M.~Combes, R.~Schrader and R.~Seiler \cite{Combes-Schrader-Seiler}.

There exists a useful interaction between capacities and the Feynman--Kac formula.
It was discussed e.g.\ in the books by I.~Chavel \cite{Chavel}, 
K.~Ito and H.~ McKean \cite{Ito-McKean}, 
M.~Kac \cite{Kac} and B.~Simon \cite{Simon79b}. 
M.~Kac and J.-M.~Luttinger \cite{Kac-Luttinger}  noticed an interesting relation between
the scattering length and capacity.

\end{document}